\documentclass{article}
\usepackage{amsfonts,amssymb, amsmath}

\begin{document}

\begin{center}
{\Large \ Some properties of symbol algebras of degree three}

\begin{equation*}
\end{equation*}%
Cristina FLAUT and Diana SAVIN%
\begin{equation*}
\end{equation*}
\end{center}

\textbf{Abstract. }{\small In this paper we study some properties of the
matrix representations of the symbol algebras of degree three. \ Moreover,
we study some equations with coefficients in these algebras, we define the
Fibonacci symbol elements and we obtain some of their properties. }

\bigskip


\begin{equation*}
\end{equation*}
\begin{equation*}
\end{equation*}

\textbf{0. Preliminaries}%
\begin{equation*}
\end{equation*}

\bigskip

Let $n$ be an arbitrary positive integer, let \ $K$ be a field whose $%
char(K) $ does not divide $n$ and contains $\omega ,~$a primitive $n$-th
root of the unity. $\ $Let $K^{\ast }=K\backslash \{0\},$ $a,b$ $\in K^{\ast
}$ and let $S$ be the algebra over $K$ generated by elements $x$ and $y$
where%
\begin{equation*}
x^{n}=a,y^{n}=b,yx=\omega xy.
\end{equation*}

This algebra is called a \textit{symbol algebra }(also known as a \textit{%
power norm residue algebra}) and it is denoted by $\left( \frac{a,~b}{%
K,\omega }\right) .$ J. Milnor, in [Mi; 71], calls it "symbol algebra"
because of its connection with the $K$-theory and with the Steinberg symbol.
For $n=2,$ we obtain the quaternion algebra. For details about Steinberg
symbol, the reader is referred to [La; 05].

In the paper [Fla; 12], using the associated trace form for a symbol
algebra, the author studied some properties of \ such objects \ and gave
some conditions for a symbol algebra to be with division or not. Starting
from these results, we intend to find some examples of division symbol
algebras. Since such an example is not easy to provide, we try to find first
sets of invertible elements in a symbol algebra and study their algebraic
properties and structure (see, for example, Proposition 4.4\textbf{)}. 

In the special cases of quaternion algebras, octonion algebras and symbol
algebras, the definition of a division algebra is equivalent to the fact
that all its nonzero elements are invertible (since these algebras admit a
type of norm $n$ that permits composition, i.e. $n(ab)=n\left( a\right)
n\left( b\right) ,$ for all elements $a,b$ in such an algebra). From this
idea, for a symbol algebra with $a=b=1$ (Theorem 4.7), we proved that all
Fibonacci symbol elements are invertible. This set is included in the set of
the elements denoted by $M$ from Proposition 4.4, in which we proved that $M$
is a $Z$-module. We intend to complete this set with other invertible
elements such that the obtained set to be an algebra not only a $Z$-module. 

The study of symbol algebras in general, and of degree three in particular,
involves very complicated calculus and, usually, can be hard to find \
examples for some notions.

From this reason, the present paper is rather technical, which is however
unavoidable seen the subject.

The paper is structurated in four sections. In the first section are
introduced \ some definitions and general properties \ of these algebras.
Since the theory of symbol algebras has many applications in different areas
of mathematics and other applied sciences, in sections two, we studied
matrix representations of symbol algebras of degree three and, in section
three, we used some of these properties to solve some equations with
coefficients in these algebras.

Since symbol algebras generalize the quaternion algebras, starting from some
results given in the paper [Ho; 63], in which the author defined \ and
studied Fibonacci quaternions, we define in a similar manner the Fibonacci
symbol elements and we study their properties. We computed the formula for
the reduced norm of a Fibonacci symbol element (Proposition 4.5) and, using
this expression, we find an infinite set of invertible elements (Theorem
4.7). We hope that this kind of sets of invertible elements will help us to
provide, in the future, examples of symbol division algebra of degree three
or of degree greater than three.

\begin{equation*}
\end{equation*}

\textbf{1.} \textbf{Introduction}%
\begin{equation*}
\end{equation*}

In the following, we assume that $K$ is a commutative field with $charK\neq
2,3$ and $A$ is a finite dimensional algebra over the field \ $K.$ The 
\textit{center} \ $C(A)$ of an algebra $A$ is the set of all elements $c\in
A $ which commute and associate with all elements $x\in $ $A.$ An algebra $A$
is a \textit{simple} algebra if \ $A$ is not a zero algebra and \ $\{0\}$
and $A$ are the only ideals of $A.$ The algebra $A$ is called \textit{%
central simple} if \ the algebra\ $A_{F}=F\otimes _{K}A$ is simple for every
field extension $F$ of $K$. Equivalently, a central simple algebra is a
simple algebra with $C(A)=K.$ We remark that each simple algebra is central
simple over its center. If $A$ is a central simple algebra, then $\dim
A=n=m^{2},$ with $m\in \mathbb{N}.$ The \textit{degree} of the central
simple algebra $A,$ denoted by $DegA,$ is $DegA=m.$

If $A$ is an algebra over the field $K\,,$ a \textit{subfield} of the
algebra $A$ is a subalgebra $L$ of $A$ such that $L$ is a field. The
subfield $L$ is called a \textit{maximal subfield} of the algebra $A$ if
there is not a subfield $F$ of $A$ such that $L\subset F.$ If the algebra $A$
is a central simple algebra, the subfield $L$ of the \ algebra $A$ is called
a\textit{\ strictly maximally subfield} of $A$ if $[L:K]=m,$ where $[L:K]$
is the degree of the extension $K\subset L.$

Let $L\subset M$ be a field extension. This extension is called a \textit{%
cyclic extension} if \ it is a Galois extension and the Galois group $%
G\left( M/L\right) $ is a cyclic group. A central simple algebra $A$ is
called a \textit{cyclic algebra } if there is $L,$ a strictly maximally
subfield of the algebra $A,$ such that $L/K$ is a cyclic extension.\medskip

\textbf{Proposition 1.1.} ([Pi; 82], Proposition a, p. 277) \textit{Let} \ $%
K\subset L$ \textit{be a cyclic extension with the Galois cyclic group} $%
G=G\left( L/K\right) $ \textit{of order} $n$ \textit{and generated \ by the
element} $\sigma .$ \textit{If} $\ A$ \textit{is a cyclic algebra and
contains} $\ L$ \textit{as a strictly maximally subfield, then there is an
element} $x\in A-\{0\}$ \textit{such that:}

\textit{i)} $A=\underset{0\leq j\leq n-1}{\bigoplus }x^{j}L;$

\textit{ii)} $x^{-1}\gamma x=\sigma \left( \gamma \right) ,$ \textit{for all}
$\gamma \in L;$

\textit{iii)} $x^{n}=a\in K^{\ast }.\medskip $

We will denote a cyclic algebra $A$ with $\left( L,\sigma ,a\right)
.\medskip $

We remark that a symbol algebra is a central simple cyclic algebra of degree 
$n.$ For details about central simple algebras and cyclic algebras, the
reader is referred to [Pi; 82].\medskip

\textbf{Definition 1.2. }[Pi; 82] Let $A$ be an algebra over the field $K.$
If $K\subset L$ is a finite field extension and $n$ a natural number, then a 
$K$-algebra morphism $\varphi :A\rightarrow \mathcal{M}_{n}\left( L\right) $
is called a \textit{representation of the algebra} $A.$ The $\varphi -$%
\textit{characteristic polynomial} of the element $a\in A$ is $P_{\varphi
}\left( X,a\right) =\det (XI_{n}-\varphi \left( a\right) ),$ the $\varphi -$%
\textit{norm} of the element $a\in A$ is $\eta _{\varphi }\left( a\right)
=\det \varphi \left( a\right) $ and the $\varphi -$\textit{trace} of the
element $a\in A$ is $\tau _{\varphi }\left( a\right) =tr\left( \varphi
\left( y\right) \right) .$ If $A$ is a $K-$ central simple algebra such that 
$n=DegA,$ then the representation $\varphi $ is called a \textit{splitting
representation} of the algebra $A.\medskip $

\textbf{Remark 1.3.} i) If $X$,$Y\in \mathcal{M}_{n}\left( K\right) ,$ $K$
an arbitrary field, then we know that $tr\left( X^{t}\right) =tr\left(
X\right) $ and $tr\left( X^{t}Y\right) =tr\left( XY^{t}\right) .$ It results
that $tr\left( XYX^{-1}\right) =tr\left( Y\right) .$

ii) ([Pi; 82], p. 296) If \ $\varphi _{1}:A\rightarrow \mathcal{M}_{n}\left(
L_{1}\right) ,~\varphi _{2}:A\rightarrow \mathcal{M}_{n}\left( L_{2}\right) $
are two splitting representations of \ the $K-$algebra $A,$ then \ $%
P_{\varphi _{1}}\left( X,a\right) =P_{\varphi _{2}}\left( X,a\right) .$ It
results that the $\varphi _{1}-$\textit{characteristic polynomial }\ is the
same with $\ $the $\varphi _{2}-$\textit{characteristic polynomial,}$\
\varphi _{1}-$norm is the same with $\varphi _{2}-$norm\textit{,}\newline
the $\varphi _{1}-$trace\textit{\ }is the same with $\varphi _{2}-$trace and
we will denote them by$\ \ P\left( X,a\right) $ instead of $\ P_{\varphi
}\left( X,a\right) ,$ $\eta _{A/K}$ instead of $\eta _{\varphi _{1}}$ or
simply $\eta $ and $\tau $ instead of $\tau _{\varphi },$ when is no
confusion in notation. In this case, the polynomial \ $P\left( X,a\right) $
is called \textit{the characteristic polynomial, }the norm $\eta $ is called
the \textit{reduced norm} of the element $a\in A$ and $\tau $ is called the 
\textit{trace} of the element $a\in A.\medskip $

\textbf{Proposition 1.4.} ([Pi; 82], Corollary a, p. 296)\textit{\ If} \ $A$ 
\textit{is a central simple algebra over the field} $K$ $\ $\ \textit{of
degree} \ $m,$ $\varphi :A\rightarrow \mathcal{M}_{r}\left( L\right) $ 
\textit{is a matrix representation of} $A$, \textit{then} $m\mid r,$ $\eta
_{\varphi }=\eta ^{r/m},$ $\tau _{\varphi }=\left( r/m\right) \tau $ \textit{%
and} \ $\eta _{\varphi }\left( a\right) ,\tau _{\varphi }\left( a\right) \in
K,$ \textit{for all} $a\in A.\medskip $

\begin{equation*}
\end{equation*}
\textbf{2. Matrix representations for the symbol algebras of degree three}%
\begin{equation*}
\end{equation*}

Let $\omega $ be a cubic root of unity, $K$ be a field such that $\omega \in
K$ and $S=\left( \frac{a,~b}{K,\omega }\right) $ be a symbol algebra over
the field $K$ generated by elements $x$ and $y$ where%
\begin{equation}
x^{3}=a,y^{3}=b,yx=\omega xy,a,b\in K^{\ast }.  \tag{2.1.}
\end{equation}

In [Ti; 00], the author gave many properties of the left and right matrix
representations for the real quaternion algebra. Since symbol algebras
generalize the quaternion algebras, using some ideas from this paper, in the
following, we will study the left and the right matrix representations for
the symbol algebras of degree $3.\medskip $

A basis in the algebra $S$ is 
\begin{equation}
B=\{1,x,x^{2},y,y^{2},xy,x^{2}y^{2},x^{2}y,xy^{2}\},  \tag{2.2.}
\end{equation}%
(see [Mil; 12], [Gi, Sz; 06 ], [Fla; 12]), [Sa; Fa; Ci; 09]).

Let $z\in S,$%
\begin{equation}
z=c_{0}+c_{1}x+c_{2}x^{2}+c_{3}y+c_{4}y^{2}+c_{5}xy+c_{6}x^{2}y^{2}+c_{7}x^{2}y+c_{8}xy^{2}
\tag{2.3.}
\end{equation}%
and $\Lambda \left( z\right) \in \mathcal{M}_{9}\left( K\right) $ be the
matrix with the coefficients in $K$ \ which its columns are the coordinates
of the elements $\{z\cdot
1,zx,zx^{2},zy,zy^{2},zxy,zx^{2}y^{2},zx^{2}y,zxy^{2}\}$\ in the basis $B:$ $%
\medskip \medskip $

$\Lambda \left( z\right) =\left( 
\begin{array}{ccccccccc}
c_{0} & ac_{2} & ac_{1} & bc_{4} & bc_{3} & ab\omega ^{2}c_{6} & ab\omega
^{2}c_{5} & ab\omega c_{8} & ab\omega c_{7} \\ 
c_{1} & c_{0} & ac_{2} & bc_{8} & bc_{5} & b\omega ^{2}c_{4} & ab\omega
^{2}c_{7} & ab\omega c_{6} & b\omega c_{3} \\ 
c_{2} & c_{1} & c_{0} & bc_{6} & bc_{7} & b\omega ^{2}c_{8} & a\omega
^{2}c_{3} & b\omega c_{4} & b\omega c_{5} \\ 
c_{3} & a\omega c_{7} & a\omega ^{2}c_{5} & c_{0} & bc_{4} & ac_{2} & 
ab\omega c_{8} & ac_{1} & ab\omega ^{2}c_{6} \\ 
c_{4} & a\omega ^{2}c_{6} & a\omega c_{8} & c_{3} & c_{0} & a\omega c_{7} & 
ac_{1} & a\omega ^{2}c_{5} & ac_{2} \\ 
c_{5} & \omega c_{3} & a\omega ^{2}c_{7} & c_{1} & bc_{8} & c_{0} & ab\omega
c_{6} & ac_{2} & b\omega ^{2}c_{4} \\ 
c_{6} & \omega ^{2}c_{8} & \omega c_{4} & c_{7} & c_{2} & \omega c_{5} & 
c_{0} & \omega ^{2}c_{3} & c_{1} \\ 
c_{7} & \omega c_{5} & \omega ^{2}c_{3} & c_{2} & bc_{6} & c_{1} & b\omega
c_{4} & c_{0} & b\omega ^{2}c_{8} \\ 
c_{8} & \omega ^{2}c_{4} & a\omega c_{6} & c_{5} & c_{1} & \omega c_{3} & 
ac_{2} & a\omega ^{2}c_{7} & c_{0}%
\end{array}%
\right) .\medskip \medskip $\newline
Let $\alpha _{ij}\in \mathcal{M}_{3}\left( K\right) $ be the matrix with $1$
in position $\left( i,j\right) $ and zero in the rest and\medskip\ \newline
$\gamma _{1}=\left( 
\begin{array}{ccc}
0 & 0 & a \\ 
1 & 0 & 0 \\ 
0 & 1 & 0%
\end{array}%
\right) ,\beta _{1}=\left( 
\begin{array}{ccc}
0 & 0 & 0 \\ 
0 & 0 & 1 \\ 
0 & 1 & 0%
\end{array}%
\right) ,\beta _{2}=\left( 
\begin{array}{ccc}
0 & 1 & 0 \\ 
1 & 0 & 0 \\ 
0 & 0 & 0%
\end{array}%
\right) ,$\newline
$\beta _{3}=\left( 
\begin{array}{ccc}
1 & 0 & 0 \\ 
0 & 0 & 0 \\ 
0 & \omega & 0%
\end{array}%
\right) ,\beta _{4}=\left( 
\begin{array}{ccc}
0 & 0 & 0 \\ 
0 & 0 & 1 \\ 
\omega & 0 & 0%
\end{array}%
\right) .\medskip \medskip $

\textbf{Proposition 2.1.} \textit{The map} $\Lambda :S\rightarrow \mathcal{M}%
_{9}\left( K\right) ,~z\mapsto \Lambda \left( z\right) $ \textit{is a} $K-$%
\textit{algebra morphism}.\medskip\ 

\textbf{Proof.} With the above notations, let $\Lambda \left( x\right)
=X=\left( 
\begin{array}{ccc}
\gamma _{1} & 0 & 0 \\ 
0 & \alpha _{31} & a\beta _{2} \\ 
0 & \beta _{1} & \alpha _{13}%
\end{array}%
\right) \in \mathcal{M}_{9}\left( K\right) $ and $\Lambda \left( y\right)
=Y=\left( 
\begin{array}{ccc}
0 & b\alpha _{12} & \omega b\beta _{4} \\ 
\beta _{3} & \alpha _{21} & 0 \\ 
\omega ^{2}\alpha _{23} & \omega \alpha _{33} & \omega ^{2}\alpha _{12}%
\end{array}%
\right) .\medskip $

By straightforward calculations, we obtain:\newline
$\Lambda \left( x^{2}\right) =X^{2},\Lambda \left( y^{2}\right)
=Y^{2},\Lambda \left( xy\right) =XY,$\newline
$\Lambda \left( x^{2}y\right) =\Lambda \left( x^{2}\right) \Lambda \left(
y\right) =\Lambda (x)\Lambda (xy)=X^{2}Y,$\newline
$\Lambda \left( xy^{2}\right) =\Lambda \left( x\right) \Lambda \left(
y^{2}\right) =\Lambda (xy)\Lambda (y)=XY^{2},$\newline
$\Lambda \left( x^{2}y^{2}\right) =\Lambda \left( x^{2}\right) \Lambda
\left( y^{2}\right) =\Lambda (x)\Lambda (xy^{2})=\Lambda (x^{2}y)\Lambda
(y)=X^{2}Y^{2}.$

Therefore, we have $\Lambda \left( z_{1}z_{2}\right) =\Lambda \left(
z_{1}\right) \Lambda \left( z_{2}\right) .$ It results that $\Lambda $ is a $%
K-$algebra morphism. $\Box \medskip $

The morphism $\Lambda $ is called the\ \textit{\ left matrix representation}
for the algebra $S.\medskip $

\textbf{Definition 2.2.} For $Z\in S,$ we denote by $\overrightarrow{Z}%
=\left( c_{0},c_{1},c_{2},c_{3},c_{4},c_{5},c_{6},c_{7},c_{8}\right) ^{t}%
\mathit{\ }\in \mathcal{M}_{9\times 1}\left( K\right) ~$\textit{the vector
representation }of the element $Z.\medskip $

\textbf{Proposition 2.3.} \textit{Let} $Z,A\in S,$ \textit{then:}

\textit{i)} $\overrightarrow{Z}=\Lambda \left( Z\right) \left( 
\begin{array}{c}
1 \\ 
0%
\end{array}%
\right) ,$ \textit{where} $0\in \mathcal{M}_{8\times 1}\left( K\right) $ 
\textit{is the zero matrix.}

\textit{ii)} $\overrightarrow{AZ}=\Lambda \left( A\right) \overrightarrow{Z}%
.\medskip $

\textbf{Proof.} \ ii) From i), we obtain that\newline
$\overrightarrow{AZ}=\Lambda \left( AZ\right) \left( 
\begin{array}{c}
1 \\ 
0%
\end{array}%
\right) =\Lambda \left( A\right) \Lambda \left( Z\right) \left( 
\begin{array}{c}
1 \\ 
0%
\end{array}%
\right) =\Lambda \left( A\right) \overrightarrow{Z}.\medskip \medskip \Box $

\textbf{Remark 2.4.} 1) We remark that an element $z\in S$ is an invertible
element in $S$ if and only if $\det \Lambda \left( z\right) \neq 0.$

2) The $\Lambda -$\textit{norm} of the element $z\in S$ is $\eta _{\Lambda
}\left( z\right) =\det \Lambda \left( z\right) =\eta ^{3}\left( z\right) $
and $\tau _{\Lambda }\left( z\right) =9tr\Lambda \left( z\right) .$ Indeed,
from Proposition 1.4, if $~A=S,K=L,m=3,r=9,\varphi =\Lambda ,$ we obtain the
above relation.

3) We have that $\tau _{\Lambda }\left( z\right) =tr\Lambda \left( z\right)
=9c_{0}.\medskip \medskip $

Let $z\in S,z=A+By+Cy^{2},$ where $%
A=c_{0}+c_{1}x+c_{2}x^{2},B=c_{3}+c_{5}x+c_{7}x^{2},C=c_{4}+c_{8}x+c_{6}x^{2}. 
$ We denote by $z_{\omega }=A+\omega By+\omega ^{2}Cy^{2},z_{\omega
^{2}}=A+\omega ^{2}By+\omega Cy^{2}.\medskip $

\textbf{Proposition 2.5. }\textit{Let} $S=\left( \frac{1,1}{K,\omega }%
\right) .$ \textit{Then} $\eta _{\Lambda }\left( z\right) =\eta _{\Lambda
}\left( z_{\omega }\right) =\eta _{\Lambda }\left( z_{\omega ^{2}}\right)
.\medskip $

\textbf{Proof. }The left matrix representation for the element $z\in S\,$\
is\medskip\ $\Lambda \left( z\right) =\left( 
\begin{array}{ccccccccc}
c_{0} & c_{2} & c_{1} & c_{4} & c_{3} & \omega ^{2}c_{6} & \omega ^{2}c_{5}
& \omega c_{8} & \omega c_{7} \\ 
c_{1} & c_{0} & c_{2} & c_{8} & c_{5} & \omega ^{2}c_{4} & \omega ^{2}c_{7}
& \omega c_{6} & \omega c_{3} \\ 
c_{2} & c_{1} & c_{0} & c_{6} & c_{7} & \omega ^{2}c_{8} & \omega ^{2}c_{3}
& \omega c_{4} & \omega c_{5} \\ 
c_{3} & \omega c_{7} & \omega ^{2}c_{5} & c_{0} & c_{4} & ac_{2} & \omega
c_{8} & c_{1} & \omega ^{2}c_{6} \\ 
c_{4} & \omega ^{2}c_{6} & \omega c_{8} & c_{3} & c_{0} & \omega c_{7} & 
c_{1} & \omega ^{2}c_{5} & c_{2} \\ 
c_{5} & \omega c_{3} & \omega ^{2}c_{7} & c_{1} & c_{8} & c_{0} & \omega
c_{6} & c_{2} & \omega ^{2}c_{4} \\ 
c_{6} & \omega ^{2}c_{8} & \omega c_{4} & c_{7} & c_{2} & \omega c_{5} & 
c_{0} & \omega ^{2}c_{3} & c_{1} \\ 
c_{7} & \omega c_{5} & \omega ^{2}c_{3} & c_{2} & c_{6} & c_{1} & \omega
c_{4} & c_{0} & \omega ^{2}c_{8} \\ 
c_{8} & \omega ^{2}c_{4} & \omega c_{6} & c_{5} & c_{1} & \omega c_{3} & 
c_{2} & \omega ^{2}c_{7} & c_{0}%
\end{array}%
\right) $ \newline
and for the element $z_{\omega }$ is\medskip \medskip \newline
\ $\Lambda \left( z_{\omega }\right) =\left( 
\begin{array}{ccccccccc}
c_{0} & c_{2} & c_{1} & \omega ^{2}c_{4} & \omega c_{3} & \omega c_{6} & 
c_{5} & c_{8} & \omega ^{2}c_{7} \\ 
c_{1} & c_{0} & c_{2} & \omega ^{2}c_{8} & \omega c_{5} & \omega c_{4} & 
c_{7} & c_{6} & \omega ^{2}c_{3} \\ 
c_{2} & c_{1} & c_{0} & \omega ^{2}c_{6} & \omega c_{7} & \omega c_{8} & 
c_{3} & c_{4} & \omega ^{2}c_{5} \\ 
c_{3} & \omega ^{2}c_{7} & c_{5} & c_{0} & \omega ^{2}c_{4} & c_{2} & c_{8}
& c_{1} & \omega c_{6} \\ 
c_{4} & \omega c_{6} & c_{8} & \omega c_{3} & c_{0} & \omega ^{2}c_{7} & 
c_{1} & c_{5} & c_{2} \\ 
c_{5} & \omega ^{2}c_{3} & c_{7} & c_{1} & \omega ^{2}c_{8} & c_{0} & c_{6}
& c_{2} & \omega c_{4} \\ 
c_{6} & \omega c_{8} & c_{4} & \omega c_{7} & c_{2} & \omega ^{2}c_{5} & 
c_{0} & c_{3} & c_{1} \\ 
c_{7} & \omega ^{2}c_{5} & c_{3} & c_{2} & \omega ^{2}c_{6} & c_{1} & c_{4}
& c_{0} & \omega c_{8} \\ 
c_{8} & \omega c_{4} & c_{6} & \omega c_{5} & c_{1} & \omega ^{2}c_{3} & 
c_{2} & c_{7} & c_{0}%
\end{array}%
\right) .\medskip $

Denoting by $D_{rs}=\left( d_{ij}^{rs}\right) \in \mathcal{M}_{9}\left(
K\right) $ the matrix defined such that $d_{kk}^{rs}=1$ for $k\notin
\{r,s\},d_{rr}^{rs}=d_{ss}^{rs}=0,$ $d_{rs}^{rs}=d_{sr}^{rs}=1$ and zero in
the rest, we have that $\det D_{rs}=-1$. If we multiply a matrix $A$ to the
left with $D_{rs},$ the new matrix is obtained from $A$ by changing the line 
$r$ with the line $s$ and if we multiply a matrix $A$ to the right with $%
D_{rs},$ the new matrix is obtained from $A$ by changing the column $r$ with
the column $s$. By straightforward calculations, it results that $\Lambda
\left( z_{\omega }\right) =D_{79}D_{48}D_{46}D_{57}D_{21}D_{23}\Lambda
\left( z\right) D_{12}D_{23}D_{49}D_{79}D_{48}D_{59},$ therefore $\eta
_{\Lambda }\left( z\right) =\det \Lambda \left( z\right) =\det \Lambda
\left( z_{\omega }\right) =\eta _{\Lambda }\left( z_{\omega }\right) .$ In
the same way, we get that $\eta _{\Lambda }\left( z\right) =\det \Lambda
\left( z\right) =\det \Lambda \left( z_{\omega ^{2}}\right) =\eta _{\Lambda
}\left( z_{\omega ^{2}}\right) .\Box \medskip $

Similar to the matrix $\Lambda \left( z\right) ,z\in S,$ we define \ $\Gamma
\left( z\right) \in \mathcal{M}_{9}\left( K\right) $ to be the matrix with
the coefficients in $K$ \ of the basis $B$ for the elements \newline
$\{z\cdot 1,xz,x^{2}z,yz,y^{2}z,xyz,x^{2}y^{2}z,x^{2}yz,xy^{2}z\}$ in their
columns$.$ This matrix is\medskip \medskip\ \newline
$\Gamma \left( z\right) =\left( 
\begin{array}{ccccccccc}
c_{0} & ac_{2} & ac_{1} & bc_{4} & bc_{3} & ab\omega ^{2}c_{6} & ab\omega
^{2}c_{5} & ab\omega c_{8} & ab\omega c_{7} \\ 
c_{1} & c_{0} & ac_{2} & b\omega c_{8} & b\omega ^{2}c_{5} & bc_{4} & 
ab\omega c_{7} & ab\omega ^{2}c_{6} & bc_{3} \\ 
c_{2} & c_{1} & c_{0} & b\omega ^{2}c_{6} & b\omega c_{7} & b\omega c_{8} & 
ac_{3} & bc_{4} & b\omega ^{2}c_{5} \\ 
c_{3} & ac_{7} & ac_{5} & c_{0} & bc_{4} & a\omega ^{2}c_{2} & ab\omega
^{2}c_{8} & a\omega c_{1} & ab\omega c_{6} \\ 
c_{4} & ac_{6} & ac_{8} & c_{3} & c_{0} & a\omega ^{2}c_{7} & a\omega
^{2}c_{1} & a\omega c_{5} & a\omega c_{2} \\ 
c_{5} & c_{3} & ac_{7} & \omega c_{1} & b\omega ^{2}c_{8} & c_{0} & ab\omega
c_{6} & a\omega ^{2}c_{2} & bc_{4} \\ 
c_{6} & c_{8} & c_{4} & \omega ^{2}c_{7} & \omega c_{2} & \omega c_{5} & 
c_{0} & c_{3} & \omega ^{2}c_{1} \\ 
c_{7} & c_{5} & c_{3} & \omega ^{2}c_{2} & b\omega c_{6} & \omega c_{1} & 
bc_{4} & c_{0} & b\omega ^{2}c_{8} \\ 
c_{8} & c_{4} & ac_{6} & \omega c_{5} & \omega ^{2}c_{1} & c_{3} & a\omega
c_{2} & a\omega ^{2}c_{7} & c_{0}%
\end{array}%
\right) \medskip \medskip $

Let $\alpha _{ij}\in \mathcal{M}_{3}\left( K\right) $ be the matrix with $1$
in position $\left( i,j\right) $ and zero in the rest and\medskip\ \newline
$\beta _{5}=\left( 
\begin{array}{ccc}
0 & 0 & 0 \\ 
0 & 0 & 1 \\ 
0 & \omega & 0%
\end{array}%
\right) ,\beta _{6}=\left( 
\begin{array}{ccc}
0 & 1 & 0 \\ 
\omega & 0 & 0 \\ 
0 & 0 & 0%
\end{array}%
\right) ,$\newline
$\beta _{7}=\left( 
\begin{array}{ccc}
1 & 0 & 0 \\ 
0 & 0 & 0 \\ 
0 & 1 & 0%
\end{array}%
\right) ,\beta _{8}=\left( 
\begin{array}{ccc}
0 & 0 & 0 \\ 
0 & 0 & 1 \\ 
1 & 0 & 0%
\end{array}%
\right) .\medskip \medskip $

\textbf{Proposition 2.6.} \textit{With the above notations, we have:}

\textit{i)} $\Gamma \left( d_{1}z_{1}+d_{2}z_{2}\right) =d_{1}\Gamma \left(
z_{1}\right) +d_{2}\Gamma \left( z_{2}\right) ,$ \textit{for all} $\
z_{1},z_{2}\in S$ \textit{and} $d_{1},d_{1}\in K.$

\textit{ii)} $\Gamma \left( z_{1}z_{2}\right) =\Gamma \left( z_{2}\right)
\Gamma \left( z_{1}\right) ,$ \textit{for all} $\ z_{1},z_{2}\in S.\medskip $

\textbf{Proof}. Let $\Gamma \left( x\right) =U=\left( 
\begin{array}{ccc}
\gamma _{1} & 0 & 0 \\ 
0 & \omega \alpha _{31} & a\omega \beta _{6} \\ 
0 & \omega \beta _{5} & \omega ^{2}\alpha _{13}%
\end{array}%
\right) \in \mathcal{M}_{9}\left( K\right) $ and \newline
$\Gamma \left( y\right) =V=\left( 
\begin{array}{ccc}
0 & b\alpha _{12} & \omega b\beta _{8} \\ 
\beta _{7} & \alpha _{21} & 0 \\ 
\alpha _{23} & \alpha _{33} & \alpha _{12}%
\end{array}%
\right) .\medskip $

By straightforward calculations, we obtain:\newline
$\Gamma \left( x^{2}\right) =U^{2},\Gamma \left( y^{2}\right) =V^{2},\Gamma
\left( yx\right) =UV,$\newline
$\Gamma \left( x^{2}y\right) =\Gamma \left( y\right) \Gamma \left(
x^{2}\right) =\Gamma (xy)\Gamma (x)=VU^{2},$\newline
$\Gamma \left( xy^{2}\right) =\Gamma \left( y^{2}\right) \Gamma \left(
x\right) =\Gamma (y)\Gamma (xy)=V^{2}U,$\newline
$\Gamma \left( x^{2}y^{2}\right) =\Gamma \left( y^{2}\right) \Gamma \left(
x^{2}\right) =\Gamma (xy^{2})\Gamma (x)=\Gamma (y)\Gamma
(x^{2}y)=V^{2}U^{2}. $

Therefore, we have that $\Gamma \left( z_{1}z_{2}\right) =\Gamma \left(
z_{2}\right) \Gamma \left( z_{1}\right) $ and $\Gamma $ is a $K-$algebra
morphism.$\medskip \Box \medskip $

\textbf{Proposition 2.7.} \textit{Let} $Z,A\in S,$ \textit{then:}

\textit{i)} $\overrightarrow{Z}=\Gamma \left( Z\right) \left( 
\begin{array}{c}
1 \\ 
0%
\end{array}%
\right) ,$ \textit{where} $0\in \mathcal{M}_{8\times 1}\left( K\right) $ 
\textit{is the zero matrix.}

\textit{ii)} $\overrightarrow{ZA}=\Gamma \left( A\right) \overrightarrow{Z}.$

iii) $\Lambda \left( A\right) \Gamma \left( B\right) =\Gamma \left( B\right)
\Lambda \left( A\right) .\medskip \Box $

\bigskip

\textbf{Proposition 2.8. }\textit{Let} $S=\left( \frac{1,1}{K,\omega }%
\right) ,$ \textit{therefore} $\eta _{\Gamma }\left( z\right) =\eta _{\Gamma
}\left( z_{\omega }\right) =\eta _{\Gamma }\left( z_{\omega ^{2}}\right)
.\medskip \medskip \Box $

\textbf{Theorem 2.9.} \ \textit{Let} $S=\left( \frac{a,b}{K,\omega }\right) $
\textit{be a symbol algebra of degree three and\newline
} $M_{9}=\left( 1,\frac{1}{a}x,\frac{1}{a}x^{2},\frac{1}{b}y,\frac{1}{b}%
y^{2},\frac{1}{ab}xy,\frac{1}{ab}x^{2}y^{2},\frac{1}{ab}x^{2}y,\frac{1}{ab}%
xy^{2}\right) ,$\ \newline
$N_{9}=\left( 1,x^{2},x,y^{2},y,x^{2}y^{2},xy,xy^{2},x^{2}y\right) ^{t}$%
\newline
$M_{10}=\left( 1,\frac{1}{a}x^{2},\frac{1}{a}x,\frac{1}{b}y^{2},\frac{1}{b}y,%
\frac{1}{ab}x^{2}y^{2},\frac{1}{ab}xy,\frac{1}{ab}xy^{2},\frac{1}{ab}%
x^{2}y\right) ,$\newline
$N_{10}=\left( 1,x,x^{2},y,y^{2},xy,x^{2}y^{2},x^{2}y,xy^{2}\right) ^{t}$. 
\textit{The following relation is true} 
\begin{equation*}
M_{9}\Lambda \left( z\right) N_{9}=M_{10}\Gamma ^{t}\left( z\right)
N_{10}=3z,z\in S.\Box
\end{equation*}%
\begin{equation*}
\end{equation*}

\textbf{3. Some equations with coefficients in a symbol algebra of degree
three} 
\begin{equation*}
\end{equation*}
Using some properties of left and right matrix representations found in the
above section, we solve some equations with coefficients in a symbol algebra
of degree three.

Let $S$ be an associative algebra of degree three. For $z\in S,$ let $%
P\left( X,z\right) $ be the characteristic polynomial for the element $a$%
\begin{equation}
P\left( X,z\right) =X^{3}-\tau \left( z\right) X^{2}+\pi \left( z\right)
X-\eta \left( z\right) \cdot 1,  \tag{3.1.}
\end{equation}%
where $\tau $ is a linear form, $~\pi $ is a quadratic form and $\eta $ a
cubic form.\medskip

\textbf{Proposition 3.1.} ([Fa; 88], Lemma) \textit{With the above
notations, denoting by} \textit{\ }$z^{\ast }=z^{2}-\tau \left( z\right)
z+\pi \left( z\right) \cdot 1,\,\,$\textit{for an associative algebra of
degree three,} \textit{we have:}

\textit{i) }$\pi \left( z\right) =\tau \left( z^{\ast }\right) .$

\textit{ii) }$2\pi \left( z\right) =\tau \left( z\right) ^{2}-\tau \left(
z^{2}\right) .$

\textit{iii) }$\tau \left( zw\right) =\tau \left( wz\right) .$

\textit{iv) }$z^{\ast \ast }=\eta \left( z\right) z.$

\textit{v) }$\left( zw\right) ^{\ast }=w^{\ast }z^{\ast }.$

\textit{vi) }$\pi \left( zw\right) =\pi \left( wz\right) .\Box \medskip
\smallskip \medskip $

In the following, we will solve some equations with coefficients in the
symbol algebra of degree three\ $S=\left( \frac{a,b}{K,\omega }\right) $.
For each element $Z\in S$ relation $\left( 3.1\right) $ holds. First, we
remark that if the element $Z\in S$ \ has $\eta \left( Z\right) \neq 0,$
then $Z$ is an invertible element. Indeed, from $\left( 3.1.\right) ,$ we
have that $ZZ^{\ast }=\eta \left( Z\right) ,$ therefore $Z^{-1}=\frac{Z}{%
\eta \left( Z\right) }.$

We consider the following equations:

\begin{equation}
AZ=ZA  \tag{3.2.}
\end{equation}%
\begin{equation}
AZ=ZB  \tag{3.3.}
\end{equation}%
\begin{equation}
AZ-ZA=C  \tag{3.4.}
\end{equation}%
\begin{equation}
AZ-ZB=C,  \tag{3.5.}
\end{equation}%
with $A,B,C\in S.\medskip $

\textbf{Proposition 3.2. }\textit{i) Equation} $\left( 3.2\right) $ \textit{%
has non-zero solutions in the algebra} $S.$

\textit{ii) If equation} $\left( 3.3\right) $ \textit{has nonzero solutions} 
$Z$ \textit{in the algebra} $S$ \textit{such that} $\eta \left( Z\right)
\neq 0$\textit{, then} $\tau \left( A\right) =\tau \left( B\right) $ \textit{%
and} $\eta \left( A\right) =\eta \left( B\right) .$

\textit{iii) If equation} $\left( 3.4\right) $ \textit{has solution, then
this solution is not unique.}

\textit{iv) If} $\Lambda \left( A\right) -\Gamma (B)$ \textit{is an
invertible matrix, then equation} $\left( 3.5\right) $ \textit{has a unique
solution.\medskip }

\textbf{Proof.} i) Using vector representation, we have that \ $%
\overrightarrow{AZ}=\overrightarrow{ZA}.$ It results that $\Lambda \left(
A\right) \overrightarrow{Z}=\Gamma \left( A\right) \overrightarrow{Z},$
therefore $(\Lambda \left( A\right) -\Gamma \left( A\right) )\overrightarrow{%
Z}=0.$ The matrix $\Lambda (A)-\Gamma \left( A\right) $ has the determinant
equal with zero (first column is zero), then equation $\left( 3.2\right) $
has non-zero solutions.

ii) We obtain $\eta \left( AZ\right) =\eta \left( BZ\right) $ and $\eta
\left( A\right) =\eta \left( B\right) $ \ since $Z$ is an invertible element
in $S.$ Using representation $\Lambda $, we get \ $\Lambda \left( A\right)
\Lambda \left( Z\right) =\Lambda \left( Z\right) \Lambda \left( B\right) .\ $%
\ From Remark 1.3 i), it results that $\Lambda \left( A\right) =\Lambda
\left( Z\right) \Lambda \left( B\right) (\Lambda \left( Z\right) )^{-1}.$
Therefore $\tau \left( A\right) =tr\left( \Lambda \left( A\right) \right)
=tr\left( \Lambda \left( Z\right) \Lambda \left( B\right) (\Lambda \left(
Z\right) )^{-1}\right) =tr\left( \Lambda \left( B\right) \right) =\tau
\left( B\right) \in K.$

iii) Using the vector representation, we have $(\Lambda \left( A\right)
-\Gamma \left( A\right) )\overrightarrow{Z}=C$ and, since the matrix $%
\Lambda (A)-\Gamma \left( A\right) $ has the determinant equal with zero
(first column is zero), if equation $\left( 3.4\right) $ has solution, then
this solution is not unique\textit{.}

iv) Using vector representation, we have $(\Lambda \left( A\right) -\Gamma
\left( B\right) )\overrightarrow{Z}=C.\Box \medskip \medskip \medskip $

\textbf{Proposition 3.3.} \ \textit{Let} \newline
$%
A=a_{0}+a_{1}x+a_{2}x^{2}+a_{3}y+a_{4}y^{2}+a_{5}xy+a_{6}x^{2}y^{2}+a_{7}x^{2}y+a_{8}xy^{2}\in S, 
$\newline
$%
B=b_{0}+b_{1}x+b_{2}x^{2}+b_{3}y+b_{4}y^{2}+b_{5}xy+b_{6}x^{2}y^{2}+b_{7}x^{2}y+b_{8}xy^{2}\in S,A_{0}=A-a_{0}\neq 0 
$ \textit{and} $B_{0}=B-b_{0}\neq 0.\ $\textit{If} $a_{0}=b_{0},A_{0}\neq
-B_{0},\eta \left( A_{0}\right) =\eta \left( B_{0}\right) =0$ \textit{and} $%
\pi \left( A_{0}\right) =\pi \left( B_{0}\right) \neq 0,$ \textit{then all
solutions of equation} $\left( 3.3\right) $ \textit{are in the }$K$\textit{%
-algebra }$\mathcal{A}\left( A,B\right) ,$ \textit{the subalgebra of }$S$ 
\textit{generated by the elements} $A$ \textit{and} $B$, \textit{and have
the form} $\lambda _{1}X_{1}+\lambda _{2}X_{2},$ \textit{where} $\lambda
_{1},\lambda _{2}\in K,$ $X_{1}=A_{0}+B_{0}$ \textit{and } $X_{2}=\pi \left(
A_{0}\right) -A_{0}B_{0}.\medskip $ \ 

\textbf{Proof.} First, we verify that $X_{1}$ and $X_{2}$ are solutions of
the equation $\left( 3.3\right) .$

Now, we prove that $X_{1}$ and $X_{2}$ are linearly independent elements. If 
$\alpha _{1}X_{1}+\alpha _{2}X_{2}=0,$ it results that $\alpha _{2}\pi
\left( A_{0}\right) =0,$ therefore $\alpha _{2}=0.$ We obtain that $\alpha
_{1}=0.$ Obviously, each element of the form $\lambda _{1}X_{1}+\lambda
_{2}X_{2},$ where $\lambda _{1},\lambda _{2}\in K$ is a solution of the
equation $\left( 3.3\right) $ and since \ $\pi \left( A_{0}\right) \neq 0$
we have that \ $\mathcal{A}\left( A,B\right) =\mathcal{A}\left(
X_{1},X_{2}\right) $. Since each solution of the equation $\left( 3.3\right) 
$ belongs to algebra $\mathcal{A}\left( A,B\right) ,$ it results that all
solutions have the form $\lambda _{1}X_{1}+\lambda _{2}X_{2},$ where $%
\lambda _{1},\lambda _{2}\in K.\Box $ 
\begin{equation*}
\end{equation*}
\bigskip \textbf{4. Fibonacci symbol elements}%
\begin{equation*}
\end{equation*}

In this section we will introduce the Fibonacci symbol elements and we will
compute the reduced norm of such an element. This relation helps us to find
an infinite set of invertible elements. First of all, we recall and give
some properties of Fibonacci numbers, properties which will be used in our
proofs.

Fibonacci numbers are the following sequence of numbers%
\begin{equation*}
0,1,1,2,3,5,8,13,21,....,
\end{equation*}%
with the $n$th term given by the formula:%
\begin{equation*}
f_{n}=f_{n-1}+f_{n-2,}\ n\geq 2,\ 
\end{equation*}%
where $f_{0}=0,f_{1}=1.$ The expression for the $n$th term is 
\begin{equation}
f_{n}=\frac{1}{\sqrt{5}}[\alpha ^{n}-\beta ^{n}],  \tag{4.1.}
\end{equation}%
where $\alpha =\frac{1+\sqrt{5}}{2}$ and $\beta =\frac{1-\sqrt{5}}{2}.$

\textbf{Remark 4.1.} \textit{Let} $(f_{n})_{n\geq 0}$ \textit{be the
Fibonacci sequence} $f_{0}=0,$ $f_{1}=1,$ $f_{n+2}=f_{n+1}+f_{n},$ $(\forall
)$ $n$$\in $ $\mathbb{N}$. \textit{Then} 
\begin{equation*}
f_{n}+f_{n+3}=2f_{n+2},(\forall )n\in \mathbb{N}.
\end{equation*}

\begin{equation*}
f_{n}+f_{n+4}=3f_{n+2},(\forall )n\in \mathbb{N}.
\end{equation*}%
\textbf{\medskip }

\textbf{Remark 4.2. \ }Let $(f_{n})_{n\geq 0}$ be the Fibonacci sequence $%
f_{0}=0,$ $f_{1}=1,$ $f_{n+2}=f_{n+1}+f_{n},$ $(\forall )$ $n$$\in $ $%
\mathbb{N}$. Let $\omega $ be a primitive root of unity of order $3,$ let $K=%
\mathbb{Q}\left( \omega \right) $ be the cyclotomic field and let $S=\left( 
\frac{a,b}{K,\omega }\right) $ be the symbol algebra of degree $3.$ Thus, $S$
has a $K$- basis $\left\{ x^{i}y^{j}|0\leq i,j<3\right\} $ such that $%
x^{3}=a $$\in $$K^{\ast },$ $y^{3}=b$$\in $$K^{\ast },$ $yx=\omega xy$ 
\newline
Let $z$$\in S$$,$ $z=\sum\limits_{i,j=1}^{n}x^{i}y^{j}c_{ij}.$ The reduced
norm of $z$ is\newline
\begin{equation*}
\eta (z)\text{=}a^{2}\cdot \left( c_{20}^{3}\text{+}bc_{21}^{3}\text{+}%
b^{2}c_{22}^{3}\text{-}3bc_{20}c_{21}c_{22}\right) \text{+}a\cdot \left(
c_{10}^{3}\text{+}bc_{11}^{3}\text{+}b^{2}c_{12}^{3}\text{-}%
3bc_{10}c_{11}c_{12}\right) \text{-}
\end{equation*}%
\begin{equation*}
\text{-}3a\cdot \left( c_{00}c_{10}c_{20}\text{+}bc_{01}c_{11}c_{21}\text{+}%
b^{2}c_{02}c_{12}c_{22}\right) \text{-}
\end{equation*}%
\begin{equation*}
-3ab\omega \left( c_{00}c_{12}c_{21}\text{+}c_{01}c_{10}c_{22}\text{+}%
c_{02}c_{11}c_{20}\right) \text{-}
\end{equation*}%
\begin{equation}
-3ab\omega ^{2}\left( c_{00}c_{11}c_{22}\text{+}c_{02}c_{10}c_{21}\text{+}%
c_{01}c_{12}c_{20}\right) \text{+}c_{00}^{3}\text{+}bc_{01}^{3}\text{+}%
b^{2}c_{02}^{3}\text{-}3bc_{00}c_{01}c_{02}.  \tag{4.2.}
\end{equation}%
(See [Pi; 82], p.299).\medskip \newline
We define the \textit{n th Fibonacci symbol element} to be the element 
\begin{equation*}
F_{n}=f_{n}\cdot 1+f_{n+1}\cdot x+f_{n+2}\cdot x^{2}+f_{n+3}\cdot y+
\end{equation*}%
\begin{equation}
+f_{n+4}\cdot xy+f_{n+5}\cdot x^{2}y+f_{n+6}\cdot y^{2}+f_{n+7}\cdot
xy^{2}+f_{n+8}\cdot x^{2}y^{2}.  \tag{4.3.}
\end{equation}%
In [Ho,61] A.F. Horadam generalized the Fibonacci numbers, giving by: 
\begin{equation*}
h_{n}=h_{n-1}+h_{n-2},\ n\in \mathbb{N},n\geq 2,
\end{equation*}%
$h_{0}=p,h_{1}=q,$ where $p,q$ are arbitrary integers. From [Ho,61],
relation 7, we have that these numbers satisfy the equality $%
h_{n+1}=pf_{n}+qf_{n+1},$ $\left( \forall \right) $$n\in \mathbb{N}$. In the
following will be easy for use to use instead of $h_{n+1}$ the notation $%
h_{n+1}^{p,q}.$ It is obviously that 
\begin{equation*}
h_{n}^{p,q}+h_{n}^{p^{^{\prime }},q^{^{\prime }}}=h_{n}^{p+p^{^{\prime
}},q+q^{^{\prime }}},
\end{equation*}%
for \ $n$ \textit{\ }a positive integers number and $p,q,p^{^{\prime
}},q^{^{\prime }}$ integers numbers. (See [Fl, Sh; 12])

\medskip

\textbf{Remark 4.3. } \textit{Let} $(f_{n})_{n\geq 0}$ \textit{be the
Fibonacci sequence} $f_{0}=0,$ $f_{1}=1,$ $f_{n+2}=f_{n+1}+f_{n},$ $(\forall
)$ $n$$\in $ $\mathbb{N}$. \textit{Then}\newline
i) $f_{n}^{2}+f_{n-1}^{2}=f_{2n-1},(\forall )n\in \mathbb{N}^{\ast };$ 
\newline
ii) $f_{n+1}^{2}-f_{n-1}^{2}=f_{2n},(\forall )n\in \mathbb{N}^{\ast };$ 
\newline
iii) $f_{n+3}^{2}=2f_{n+2}^{2}+2f_{n+1}^{2}-f_{n}^{2},(\forall )n\in \mathbb{%
N};$ \newline
iv) $f_{n}^{2}-f_{n-1}f_{n+1}=\left( -1\right) ^{n-1},(\forall )n\in \mathbb{%
N}^{\ast };$ \newline
v) $f_{2n}=f_{n}^{2}+2f_{n}f_{n-1}$\newline
\bigskip

We define the $n$\textit{th generalized Fibonacci symbol element} to be the
element 
\begin{equation*}
H_{n}^{p,q}=h_{n}^{p,q}\cdot 1+h_{n+1}^{p,q}\cdot x+h_{n+2}^{p,q}\cdot
x^{2}+h_{n+3}^{p,q}\cdot y+
\end{equation*}%
\begin{equation}
+h_{n+4}^{p,q}\cdot xy+h_{n+5}^{p,q}\cdot x^{2}y+h_{n+6}^{p,q}\cdot
y^{2}+h_{n+7}^{p,q}\cdot xy^{2}+h_{n+8}^{p,q}\cdot x^{2}y^{2}.  \tag{4.4.}
\end{equation}%
What algebraic structure has the set of Fibonacci symbol elements or the set
of generalized Fibonacci symbol? The answer will be found in the following
proposition.\medskip\ \newline
\smallskip \newline
\textbf{Proposition 4.4. } \textit{Let} $M$ \textit{the set} 
\begin{equation*}
M=\left\{ \sum\limits_{i=1}^{n}H_{n_{i}}^{p_{i},q_{i}}|n\in \mathbb{N}^{\ast
},p_{i},q_{i}\in \mathbb{Z},(\forall )i=\overline{1,n}\right\} \cup \left\{
0\right\}
\end{equation*}%
\textit{is a} $\mathbb{Z}-$ \textit{module.}\newline
\smallskip \newline
\textbf{Proof.} The proof is immediate if we remark that for $n_{1},n_{2}$$%
\in $$\mathbb{N},$ $p_{1},p_{2},\alpha _{1},\alpha _{2}$$\in $$\mathbb{Z},$
we have: 
\begin{equation*}
\alpha _{1}h_{n_{1}}^{p_{1},q_{1}}+\alpha
_{2}h_{n_{2}}^{p_{2},q_{2}}=h_{n_{1}}^{\alpha _{1}p_{1},\alpha
_{1}q_{1}}+h_{n_{2}}^{\alpha _{2}p_{2},\alpha _{2}q_{2}}.
\end{equation*}%
Therefore, we obtain that $M$ is a $\mathbb{Z}-$ submodule of the symbol
algebra $S.\Box \medskip \medskip \medskip $ \smallskip \newline
In the following we will compute the reduced norm for the $n$th Fibonacci
symbol element. \newline
\smallskip \newline
\textbf{Proposition 4.5.} \textit{Let} $F_{n}$ \textit{be the }$n$\textit{th
Fibonacci symbol element. Let }$\omega $ \textit{be a primitive root of
order }$3$\textit{\ of unity and let }$K=\mathbb{Q}\left( \omega \right) $ 
\textit{be} \textit{the cyclotomic field.\ Then the norm of} $F_{n}$ \textit{%
is} \newline
\begin{equation*}
\eta (F_{n})=4a^{2}h_{n+3}^{211,14}\cdot \left(
h_{2n}^{84,135}-2f_{n}^{2}\right) +8h_{n+3}^{8,3}\cdot \left(
h_{2n}^{12,20}-f_{n}^{2}\right) +
\end{equation*}%
\begin{equation*}
\text{+}a[f_{n+2}\left( \omega h_{2n}^{30766,27923}\text{+}%
h_{2n}^{22358,20533}\right) \text{+}f_{n+3}\cdot \left( \omega
h_{2n}^{4368,1453}\text{+}h_{2n}^{14128,12163}\right) \text{-}
\end{equation*}%
\begin{equation}
-f_{n}^{2}\left( \omega h_{n+3}^{45013,22563}\text{+}h_{n+3}^{33683,27523}%
\right) -\left( -1\right) ^{n}\left( \omega h_{n+3}^{1472,26448}\text{+}%
h_{n+3}^{12982,24138}\right) ].  \tag{4.5.}
\end{equation}%
\smallskip \newline
\textbf{Proof.} In this proof, we denote with $E\left( x,y,z\right)
=x^{3}+y^{3}+z^{3}-3xyz.$ We obtain: 
\begin{equation*}
\eta (F_{n})=\eta (F_{n})+3\cdot \left( 1+\omega +\omega ^{2}\right) \cdot
\left( f_{n}^{3}+f_{n+1}^{3}+f_{n+2}^{3}\text{+...+}f_{n+8}^{3}\right) =
\end{equation*}%
\begin{equation*}
a^{2}\cdot \left( f_{n+2}^{3}\text{+}f_{n+5}^{3}\text{+}f_{n+8}^{3}\text{-}%
3f_{n+2}f_{n+5}f_{n+8}\right) \text{+}a\cdot \left( f_{n+1}^{3}\text{+}%
f_{n+4}^{3}\text{+}f_{n+7}^{3}\text{-}3f_{n+1}f_{n+3}f_{n+7}\right) \text{+}
\end{equation*}%
\begin{equation*}
\text{+}a\cdot \left( f_{n}^{3}\text{+}f_{n+1}^{3}\text{+}f_{n+2}^{3}\text{-}%
3f_{n}f_{n+1}f_{n+2}\right) \text{+}a\cdot \left( f_{n+3}^{3}\text{+}%
f_{n+4}^{3}\text{+}f_{n+5}^{3}\text{-}3f_{n+3}f_{n+4}f_{n+5}\right) \text{+}
\end{equation*}%
\begin{equation*}
\text{+}a\cdot \left( f_{n+6}^{3}\text{+}f_{n+7}^{3}\text{+}f_{n+8}^{3}\text{%
-}3f_{n+6}f_{n+7}f_{n+8}\right) \text{+}a\cdot \left( \omega ^{3}f_{n}^{3}%
\text{+}f_{n+5}^{3}+f_{n+7}^{3}\text{-}3\omega f_{n}f_{n+5}f_{n+7}\right) 
\text{+}
\end{equation*}%
\begin{equation*}
\text{+}a\cdot \left( f_{n+1}^{3}\text{+}f_{n+3}^{3}\text{+}\omega
^{3}f_{n+8}^{3}\text{-}3\omega f_{n+1}f_{n+3}f_{n+8}\right) \text{+}a\cdot
\left( f_{n+2}^{3}\text{+}f_{n+4}^{3}\text{+}\omega ^{3}f_{n+6}^{3}\text{-}%
3\omega f_{n+2}f_{n+4}f_{n+6}\right)
\end{equation*}%
\begin{equation*}
\text{+}a\cdot \left( f_{n}^{3}\text{+}f_{n+4}^{3}\text{+}\omega
^{6}f_{n+8}^{3}\text{-}3\omega ^{2}f_{n}f_{n+4}f_{n+8}\right) \text{+}a\cdot
\left( \omega ^{6}f_{n+1}^{3}\text{+}f_{n+5}^{3}\text{+}f_{n+6}^{3}\text{-}%
3\omega ^{2}f_{n+1}f_{n+5}f_{n+6}\right)
\end{equation*}%
\begin{equation*}
+a\cdot \left( f_{n+2}^{3}\text{+}f_{n+3}^{3}\text{+}\omega ^{6}f_{n+7}^{3}%
\text{-}3\omega ^{2}f_{n+2}f_{n+3}f_{n+7}\right) \text{+}\left( f_{n}^{3}%
\text{+}f_{n+3}^{3}\text{+}f_{n+6}^{3}\text{-}3f_{n}f_{n+3}f_{n+6}\right) 
\text{=}
\end{equation*}%
\begin{equation*}
=a^{2}\cdot E\left( f_{n+2},f_{n+5},f_{n+8}\right) +a\cdot E\left(
f_{n+1},f_{n+4},f_{n+7}\right) +a\cdot E\left( f_{n},f_{n+1},f_{n+2}\right) +
\end{equation*}%
\begin{equation*}
+a\cdot E\left( f_{n+3},f_{n+4},f_{n+5}\right) +a\cdot E\left(
f_{n+6},f_{n+7},f_{n+8}\right) +a\cdot E\left( \omega
f_{n},f_{n+5},f_{n+7}\right) +
\end{equation*}%
\begin{equation*}
\text{+}a\cdot E\left( f_{n+1},f_{n+3},\omega f_{n+8}\right) \text{+}a\cdot
E\left( f_{n+2},f_{n+4},\omega f_{n+6}\right) \text{+}a\cdot E\left(
f_{n},f_{n+4},\omega ^{2}f_{n+8}\right) +
\end{equation*}%
\begin{equation*}
+a\cdot E\left( \omega ^{2}f_{n+1},f_{n+5},f_{n+6}\right) +a\cdot E\left(
f_{n+2},f_{n+3},\omega ^{2}f_{n+7}\right) +E\left(
f_{n},f_{n+3},f_{n+6}\right) .
\end{equation*}%
Now, we compute $E\left( f_{n+2},f_{n+5},f_{n+8}\right) .$\newline
\begin{equation*}
E\left( f_{n+2},f_{n+5},f_{n+8}\right) =
\end{equation*}%
\begin{equation*}
\text{=}\frac{1}{2}\left( f_{n+2}\text{+}f_{n+5}\text{+}f_{n+8}\right) \left[
\left( f_{n+5}\text{-}f_{n+2}\right) ^{2}\text{+}\left( f_{n+8}\text{-}%
f_{n+5}\right) ^{2}\text{+}\left( f_{n+8}\text{-}f_{n+2}\right) ^{2}\right] .
\end{equation*}%
Using Remark 4.1, Remark 4.2, Remark 4.3 (iii) and the recurrence of the
Fibonacci sequence, we obtain: 
\begin{equation*}
E\left( f_{n+2},f_{n+5},f_{n+8}\right) =
\end{equation*}%
\begin{equation*}
\text{=}\frac{1}{2}\left( 2f_{n+4}\text{+}f_{n+8}\right) \left[ \left(
f_{n+4}\text{+}f_{n+1}\right) ^{2}\text{+}\left( f_{n+7}\text{+}%
f_{n+4}\right) ^{2}\text{+}\left( f_{n+1}\text{+}2f_{n+4}\text{+}%
f_{n+7}\right) ^{2}\right] \text{=}
\end{equation*}%
\begin{equation*}
=\frac{1}{2}\left( f_{n+4}+3f_{n+6}\right) \cdot \left[ \left(
2f_{n+3}\right) ^{2}+\left( 2f_{n+6}\right) ^{2}+\left(
2f_{n+3}+2f_{n+6}\right) ^{2}\right] =
\end{equation*}%
\begin{equation*}
=4\left( 11f_{n+2}+14f_{n+3}\right) \cdot \left(
135f_{n+1}^{2}+82f_{n}^{2}-51f_{n-1}^{2}\right) .
\end{equation*}%
Then, we have:\newline
\begin{equation}
E\left( f_{n+2},f_{n+5},f_{n+8}\right) \text{=}4\left( 11f_{n+2}\text{+}%
14f_{n+3}\right) \cdot \left( 135f_{n+1}^{2}\text{+}82f_{n}^{2}\text{-}%
51f_{n-1}^{2}\right) .  \tag{4.6.}
\end{equation}%
Using Remark 4.3 (i,ii) and the definition of the generalized Fibonacci
sequence it is easy to compute that: 
\begin{equation*}
E\left( f_{n+2},f_{n+5},f_{n+8}\right) \text{=}4\left( 11f_{n+2}\text{+}%
14f_{n+3}\right) \cdot \left( 135f_{2n}\text{+}84f_{2n-1}\text{-}%
2f_{n}^{2}\right) =
\end{equation*}%
\begin{equation*}
\text{=}4\left( 11f_{n+2}+14f_{n+3}\right) \cdot \left(
h_{2n}^{84,135}-2f_{n}^{2}\right) \text{=}4h_{n+3}^{11,14}\cdot \left(
h_{2n}^{84,135}-2f_{n}^{2}\right) .
\end{equation*}

Therefore, we obtain 
\begin{equation}
E\left( f_{n+2},f_{n+5},f_{n+8}\right) =4h_{n+3}^{11,14}\cdot \left(
h_{2n}^{84,135}-2f_{n}^{2}\right) .  \tag{4.7.}
\end{equation}

Replacing $n\rightarrow n-1$ in relation (4.6) and using Remark 4.3, (iii)
and the recurrence of the Fibonacci sequence, we obtain: 
\begin{equation}
E\left( f_{n+1},f_{n+4},f_{n+7}\right) \text{=}4\left( 3f_{n+2}\text{+}%
11f_{n+3}\right) \cdot \left( 51f_{n+1}^{2}\text{+}33f_{n}^{2}\text{-}%
20f_{n-1}^{2}\right) .\   \tag{4.8.}
\end{equation}%
Now, we calculate $E\left( f_{n},f_{n+1},f_{n+2}\right) .$ 
\begin{equation*}
E\left( f_{n},f_{n+1},f_{n+2}\right) =\frac{1}{2}\left( f_{n}\text{+}f_{n+1}%
\text{+}f_{n+2}\right) \cdot \left[ \left( f_{n+1}\text{-}f_{n}\right) ^{2}%
\text{+}\left( f_{n+2}\text{-}f_{n+1}\right) ^{2}\text{+}\left( f_{n+2}\text{%
-}f_{n}\right) ^{2}\right] \text{=}
\end{equation*}%
\begin{equation*}
\text{=}f_{n+2}\cdot \left( f_{n-1}^{2}+f_{n}^{2}+f_{n+1}^{2}\right) .
\end{equation*}%
So, we obtain $E\left( f_{n},f_{n+1},f_{n+2}\right) =f_{n+2}\cdot \left(
f_{n+1}^{2}+f_{n}^{2}+f_{n-1}^{2}\right) $ . \ \ \ \ \ \ \ \ \ \ \ (4.9.)%
\newline
Replacing $n\rightarrow n+3$ in relation (4.9.), using Remark 4.3 (iii) and
the recurrence of the Fibonacci sequence we obtain: 
\begin{equation}
E\left( f_{n+3},f_{n+4},f_{n+5}\right) =\left( f_{n+2}+2f_{n+3}\right) \cdot
\left( 23f_{n+1}^{2}+15f_{n}^{2}-9f_{n-1}^{2}\right) .\   \tag{4.10.}
\end{equation}%
Replacing $n\rightarrow n+3$ in relation (4.10), using Remark 4.3 (iii) and
the recurrence of the Fibonacci sequence we obtain: 
\begin{equation}
E\left( f_{n+6},f_{n+7},f_{n+8}\right) \text{=}\left( 3f_{n+2}\text{+}%
4f_{n+3}\right) \cdot \left( 635f_{n+1}^{2}\text{+}387f_{n}^{2}\text{-}%
239f_{n-1}^{2}\right) .\   \tag{4.11.}
\end{equation}%
Adding equalities (4.8),(4.9),(4.10),(4.11), after straightforward
calculation, we have: 
\begin{equation*}
E\left( f_{n+1},f_{n+4},f_{n+7}\right) \text{+}E\left(
f_{n},f_{n+1},f_{n+2}\right) \text{+}E\left( f_{n+3},f_{n+4},f_{n+5}\right) 
\text{+}E\left( f_{n+6},f_{n+7},f_{n+8}\right)
\end{equation*}%
\begin{equation*}
\text{=}f_{n+2}\cdot \left( 2511f_{n+1}^{2}\text{+}1573f_{n}^{2}\text{-}%
965f_{n-1}^{2}\right) \text{+}f_{n+3}\cdot \left( 4790f_{n+1}^{2}\text{+}%
3030f_{n}^{2}\text{-}1854f_{n-1}^{2}\right) .
\end{equation*}%
Using Remark 4.3 (i, ii), it is easy to compute that: 
\begin{equation*}
E\left( f_{n+1},f_{n+4},f_{n+7}\right) \text{+}E\left(
f_{n},f_{n+1},f_{n+2}\right) \text{+}E\left( f_{n+3},f_{n+4},f_{n+5}\right) 
\text{+}E\left( f_{n+6},f_{n+7},f_{n+8}\right) \text{=}
\end{equation*}%
\begin{equation*}
f_{n+2}\cdot \left( 1546f_{2n+1}\text{+}965f_{2n}\text{+}27f_{n}^{2}\right) 
\text{+}f_{n+3}\cdot \left( 2936f_{2n+1}\text{+}1854f_{2n-1}\text{+}%
94f_{n}^{2}\right) \text{=}
\end{equation*}%
\begin{equation}
=f_{n+2}\cdot h_{2n+1}^{965,1546}+f_{n+3}\cdot
h_{2n+1}^{1854,2936}+27f_{n}^{2}\cdot h_{n+4}^{67,1}.  \tag{4.12.}
\end{equation}%
Then, we obtained that: 
\begin{equation*}
E\left( f_{n+1},f_{n+4},f_{n+7}\right) \text{+}E\left(
f_{n},f_{n+1},f_{n+2}\right) \text{+}E\left( f_{n+3},f_{n+4},f_{n+5}\right) 
\text{+}E\left( f_{n+6},f_{n+7},f_{n+8}\right) \text{=}
\end{equation*}%
\begin{equation}
=f_{n+2}\cdot h_{2n+1}^{965,1546}\text{+}f_{n+3}\cdot h_{2n+1}^{1854,2936}%
\text{+}27f_{n}^{2}\cdot h_{n+4}^{67,1}.\   \tag{4.13.}
\end{equation}%
Replacing $n\ $by $n-1$ in relation (4.8), using Remark 4.3 (iii) and the
recurrence of the Fibonacci sequence we obtain: 
\begin{equation}
E\left( f_{n},f_{n+3},f_{n+6}\right) \text{=}8\left( 8f_{n+2}\text{+}%
3f_{n+3}\right) \cdot \left( 20f_{n+1}^{2}\text{+}11f_{n}^{2}-8f_{n-1}^{2}%
\right) .\   \tag{4.14.}
\end{equation}%
Using Remark 4.3 (i,ii), the definition of the generalized Fibonacci
sequence, we have: 
\begin{equation*}
E\left( f_{n},f_{n+3},f_{n+6}\right) \text{=}8\left( 8f_{n+2}\text{+}%
3f_{n+3}\right) \cdot \left[ 20\left( f_{n+1}^{2}-f_{n-1}^{2}\right) \text{+}%
12\left( f_{n-1}^{2}\text{+}f_{n}^{2}\right) -f_{n}^{2}\right] \text{=}
\end{equation*}%
\begin{equation*}
8\left( 8f_{n+2}\text{+}3f_{n+3}\right) \cdot \left( 20f_{2n}\text{+}%
12f_{2n-1}-f_{n}^{2}\right) \text{=}8h_{n+3}^{8,3}\cdot \left(
h_{2n}^{12,20}-f_{n}^{2}\right) .
\end{equation*}%
Therefore 
\begin{equation}
E\left( f_{n},f_{n+3},f_{n+6}\right) =8h_{n+3}^{8,3}\cdot \left(
h_{2n}^{12,20}-f_{n}^{2}\right) .\ \newline
\tag{4.15}
\end{equation}%
Now, we compute $E\left( \omega f_{n},f_{n+5},f_{n+7}\right) .$\newline
$E\left( \omega f_{n},f_{n+5},f_{n+7}\right) $=$\frac{1}{2}\left( \omega
f_{n}\text{+}f_{n+5}\text{+}f_{n+7}\right) \cdot \left[ \left( f_{n+7}\text{-%
}f_{n+5}\right) ^{2}\text{+}\left( f_{n+5}\text{-}\omega f_{n}\right) ^{2}%
\text{+}\left( f_{n+7}\text{-}\omega f_{n}\right) ^{2}\right] .$ By
repeatedly using of Fibonacci sequence recurrence and Remark 4.3 (iii), we
obtain: 
\begin{equation*}
E\left( \omega f_{n},f_{n+5},f_{n+7}\right) \text{=}\frac{1}{2}\left[
2\left( 2\text{+}\omega \right) f_{n+2}\text{+}\left( 7-\omega \right)
f_{n+3}\right] \cdot
\end{equation*}%
\begin{equation*}
\cdot \lbrack 104f_{n+1}^{2}\text{+}65f_{n}^{2}-40f_{n-1}^{2}\text{+}\left(
(8-\omega )f_{n+1}\text{+}(\omega -3)f_{n-1}\right) ^{2}\text{+}\left(
(11-\omega )f_{n+1}\text{+}(\omega -8)f_{n-1}\right) ^{2}].
\end{equation*}%
In the same way, we have: 
\begin{equation*}
E\left( \omega f_{n},f_{n+5},f_{n+7}\right) =\frac{1}{2}\left[ 2\left( 2%
\text{+}\omega \right) f_{n+2}+\left( 7-\omega \right) f_{n+3}\right] \cdot
\end{equation*}%
\begin{equation}
\cdot \left[ \left( \text{-}40\omega \text{+}287\right) f_{n+1}^{2}\text{+}%
\left( -24\omega \text{+}31\right) f_{n-1}^{2}\text{+}\left( -64\omega \text{%
+}285\right) f_{n}^{2}\text{-}4\left( 16\omega \text{-}55\right) \cdot
\left( \text{-}1\right) ^{n}\right] .\   \tag{4.16.}
\end{equation}%
Now, we calculate $E\left( f_{n+1},f_{n+3},\omega f_{n+8}\right) .$\newline
\begin{equation*}
E\left( f_{n+1},f_{n+3},\omega f_{n+8}\right) =\frac{1}{2}\left(
f_{n+1}+f_{n+3}+\omega f_{n+8}\right) \cdot
\end{equation*}%
\begin{equation*}
\cdot \left[ \left( f_{n+3}-f_{n+1}\right) ^{2}+\left( \omega
f_{n+8}-f_{n+3}\right) ^{2}+\left( \omega f_{n+8}-f_{n+1}\right) ^{2}\right]
.
\end{equation*}%
By repeatedly using of Fibonacci sequence recurrence and Remark 4.3 (iii),
it is easy to compute that: 
\begin{equation*}
E\left( f_{n+1},f_{n+3},\omega f_{n+8}\right) =\left[ \left( 5\omega
-1\right) f_{n+2}+2\left( 4\omega -1\right) f_{n+3}\right] \cdot
\end{equation*}%
\begin{equation}
\lbrack \text{-}2\left( 696\omega \text{+}578\right) f_{n+1}^{2}\text{-}%
\left( 182\omega \text{+}169\right) f_{n-1}^{2}\text{+}2\left( 485\omega 
\text{+}441\right) f_{n}^{2}\text{+}\left( 970\omega \text{+}881\right)
\cdot \left( \text{-}1\right) ^{n}].\   \tag{4.17.}
\end{equation}%
Now, we calculate $E\left( f_{n+2},f_{n+4},\omega f_{n+6}\right) .$ 
\begin{equation*}
E\left( f_{n+2},f_{n+4},\omega f_{n+6}\right) =
\end{equation*}%
\begin{equation*}
\text{=}\frac{1}{2}\left( f_{n+2}\text{+}f_{n+4}\text{+}\omega
f_{n+6}\right) \cdot \left[ \left( f_{n+4}\text{-}f_{n+2}\right) ^{2}\text{+}%
\left( \omega f_{n+6}\text{-}f_{n+4}\right) ^{2}\text{+}\left( \omega f_{n+6}%
\text{-}f_{n+2}\right) ^{2}\right] .
\end{equation*}%
In the same way, we obtain: 
\begin{equation*}
E\left( f_{n+2},f_{n+4},\omega f_{n+6}\right) \text{=}\frac{-1}{2}\left[
2\left( \omega \text{+}1\right) f_{n+2}\text{+}\left( 3\omega \text{+}%
1\right) f_{n+3}\right] \cdot
\end{equation*}%
\begin{equation}
\cdot \lbrack \left( 520\omega \text{+}309\right) f_{n+1}^{2}\text{+}\left(
80\omega \text{+}47\right) f_{n-1}^{2}\text{+}\left( 112\omega -21\right)
f_{n}^{2}\text{+}8\left( 14\omega \text{+}3\right) \cdot \left( -1\right)
^{n}].\   \tag{4.18.}
\end{equation}%
Adding relations (4.16), (4.17), (4.18) and making some calculus, it
results: 
\begin{equation*}
E\left( \omega f_{n},f_{n+5},f_{n+7}\right) +E\left( f_{n+1},f_{n+3},\omega
f_{n+8}\right) +E\left( f_{n+2},f_{n+4},\omega f_{n+6}\right) =
\end{equation*}%
\begin{equation*}
=f_{n+2}[\left( 10138\omega \text{+}6127\right) f_{n+1}^{2}\text{-}\left(
1210\omega \text{+}5231\right) f_{n}^{2}\text{+}\left( 301\omega \text{+}%
1132\right) f_{n-1}^{2}-
\end{equation*}%
\begin{equation*}
\text{-}\left( \text{-}1\right) ^{n}\left( 1235\omega \text{+}5315\right) ]%
\text{+}f_{n+3}[\left( 4103\omega \text{+}13544\right) f_{n+1}^{2}\text{-}%
\left( 3148\omega \text{+}8566\right) f_{n}^{2}\text{+}
\end{equation*}%
\begin{equation}
+\left( 307\omega \text{+}1771\right) f_{n-1}^{2}-\left( -1\right)
^{n}\left( 11396\omega \text{+}16279\right) ]\   \tag{4.19.}
\end{equation}%
Now, we compute $E\left( f_{n},f_{n+4},\omega ^{2}f_{n+8}\right) .$ 
\begin{equation*}
E\left( f_{n},f_{n+4},\omega ^{2}f_{n+8}\right) =
\end{equation*}%
\begin{equation*}
=\frac{1}{2}\left( f_{n}+f_{n+4}\text{+}\omega ^{2}f_{n+8}\right) \cdot %
\left[ \left( f_{n+4}\text{-}f_{n}\right) ^{2}\text{+}\left( \omega
^{2}f_{n+8}\text{-}f_{n+4}\right) ^{2}\text{+}\left( \omega ^{2}f_{n+8}\text{%
-}f_{n}\right) ^{2}\right] .
\end{equation*}%
Using Remark 4.2, the recurrence of the Fibonacci sequence we obtain: 
\begin{equation*}
f_{n}\text{+}f_{n+4}\text{=}3f_{n+2};f_{n+4}-f_{n}=3f_{n+1}\text{+}%
f_{n};f_{n+8}=5f_{n+2}\text{+}8f_{n+3}.
\end{equation*}%
Using Remark 4.3. (iv) it is easy to compute that: 
\begin{equation*}
E\left( f_{n},f_{n+4},\omega ^{2}f_{n+8}\right) \text{=}\left[ \left(
-5\omega +8\right) f_{n+2}\text{+}8\left( -\omega \text{+}1\right) f_{n+3}%
\right] \cdot \lbrack \left( 1460\omega \text{+}325\right) f_{n+1}^{2}\text{+%
}
\end{equation*}%
\begin{equation}
\text{+}\left( 208\omega \text{+}42\right) f_{n-1}^{2}-\left( 1030\omega 
\text{+}127\right) f_{n}^{2}\text{-}\left( 1030\omega \text{+}127\right)
\cdot \left( \text{-}1\right) ^{n}].\   \tag{4.20.}
\end{equation}%
Similarly it is easy to compute that: 
\begin{equation*}
E\left( f_{n+2},f_{n+3},\omega ^{2}f_{n+7}\right) =-\left[ \left( 3\omega
+2\right) f_{n+2}+\left( 5\omega +4\right) f_{n+3}\right] \cdot
\end{equation*}%
\begin{equation*}
\cdot \lbrack \left( 546\omega \text{+}112\right) f_{n+1}^{2}\text{+}\left(
80\omega \text{+}17\right) f_{n-1}^{2}-\left( 418\omega \text{+}87\right)
f_{n}^{2}\text{-}\left( 418\omega \text{+}87\right) \cdot \left( -1\right)
^{n}].
\end{equation*}%
Now, we compute $E\left( f_{n+5},f_{n+6},\omega ^{2}f_{n+1}\right) .$ 
\begin{equation*}
E\left( f_{n+5},f_{n+6},\omega ^{2}f_{n+1}\right) =\frac{1}{2}\left(
f_{n+5}+f_{n+6}+\omega ^{2}f_{n+1}\right) \cdot
\end{equation*}%
\begin{equation*}
\cdot \left[ \left( f_{n+6}-f_{n+5}\right) ^{2}+\left( f_{n+5}-\omega
^{2}f_{n+1}\right) ^{2}+\left( f_{n+6}-\omega ^{2}f_{n+1}\right) ^{2}\right]
.
\end{equation*}%
Using the recurrence of the Fibonacci sequence, we have: 
\begin{equation*}
f_{n+5}\text{+}f_{n+6}\text{+}\omega ^{2}f_{n+1}\text{=}\left( \omega \text{+%
}4\right) f_{n+2}\text{+}\left( \text{-}\omega \text{+}4\right)
f_{n+3};f_{n+6}\text{-}f_{n+5}\text{=}3f_{n+1}\text{+}2f_{n},
\end{equation*}%
\begin{equation*}
f_{n+5}\text{-}\omega ^{2}f_{n+1}\text{=}\left( \omega \text{+}6\right)
f_{n+1}\text{+}3f_{n};f_{n+6}\text{-}\omega ^{2}f_{n+1}\text{=}\left( \omega 
\text{+}9\right) f_{n+1}\text{+}4f_{n}.
\end{equation*}%
Using Remark 4.3 (iv) it is easy to compute that: 
\begin{equation*}
E\left( f_{n+5},f_{n+6},\omega ^{2}f_{n+1}\right) =\left[ \left( \omega
+4\right) f_{n+2}+\left( -\omega +4\right) f_{n+3}\right] \cdot
\end{equation*}%
\begin{equation}
\cdot \left[ \left( 23\omega \text{+}151\right) f_{n+1}^{2}\text{+}%
19f_{n-1}^{2}\text{-}\left( 6\omega \text{+}107\right) f_{n}^{2}\text{-}%
\left( 6\omega \text{+}107\right) \cdot \left( \text{-}1\right) ^{n}\right]
.\   \tag{4.22.}
\end{equation}%
Adding equalities $(4.20),(4.21),(4.22),$ we have: 
\begin{equation*}
E\left( f_{n},f_{n+4},\omega ^{2}f_{n+8}\right) \text{+}E\left( \omega
^{2}f_{n+1},f_{n+5},f_{n+6}\right) \text{+}E\left( f_{n+2},f_{n+3},\omega
^{2}f_{n+7}\right) \text{=}
\end{equation*}%
\begin{equation*}
=f_{n+2}\cdot \lbrack \left( 17785\omega \text{+}11895\right) f_{n+1}^{2}%
\text{-}\left( 13037\omega \text{+}7667\right) f_{n}^{2}\text{+}\left(
2542\omega \text{+}1658\right) f_{n-1}^{2}\text{-}
\end{equation*}%
\begin{equation*}
\text{-}\left( 13037\omega \text{+}7667\right) \cdot \left( \text{-}1\right)
^{n}]\text{+}f_{n+3}\cdot \lbrack \left( \text{-}2650\omega \text{-}%
6171\right) f_{n+1}^{2}\text{-}\left( 15052\omega \text{+}7859\right)
f_{n}^{2}
\end{equation*}%
\begin{equation}
+\left( 2608\omega +2048\right) f_{n-1}^{2}-\left( 15052\omega +7859\right)
\cdot \left( -1\right) ^{n}].  \tag{4.23.}
\end{equation}%
Adding (4.12) (4.19),(4.23), \ it results: 
\begin{equation*}
E\left( f_{n+1},f_{n+4},f_{n+7}\right) \text{+}E\left(
f_{n},f_{n+1},f_{n+2}\right) \text{+}E\left( f_{n+3},f_{n+4},f_{n+5}\right) 
\text{+}E\left( f_{n+6},f_{n+7},f_{n+8}\right)
\end{equation*}%
\begin{equation*}
\text{+}E\left( \omega f_{n},f_{n+5},f_{n+7}\right) +E\left(
f_{n+1},f_{n+3},\omega f_{n+8}\right) +E\left( f_{n+2},f_{n+4},\omega
f_{n+6}\right) +
\end{equation*}%
\begin{equation*}
E\left( f_{n},f_{n+4},\omega ^{2}f_{n+8}\right) +E\left( \omega
^{2}f_{n+1},f_{n+5},f_{n+6}\right) +E\left( f_{n+2},f_{n+3},\omega
^{2}f_{n+7}\right) \text{=}
\end{equation*}%
\begin{equation*}
f_{n+2}\cdot \left( 2511f_{n+1}^{2}\text{+}1573f_{n}^{2}-965f_{n-1}^{2}%
\right) \text{+}f_{n+3}\cdot \left( 4790f_{n+1}^{2}\text{+}3030f_{n}^{2}%
\text{-}1854f_{n-1}^{2}\right) \text{+}
\end{equation*}%
\begin{equation*}
f_{n+2}\cdot \lbrack \left( 27923\omega \text{+}18022\right) f_{n+1}^{2}%
\text{-}\left( 14247\omega \text{+}12898\right) f_{n}^{2}\text{+}\left(
2843\omega \text{+}2790\right) f_{n-1}^{2}\text{-}
\end{equation*}%
\begin{equation*}
\text{-}\left( 14272\omega \text{+}12928\right) \cdot \left( \text{-}%
1\right) ^{n}]\text{+}f_{n+3}\cdot \lbrack \left( 1453\omega \text{+}%
7373\right) f_{n+1}^{2}\text{-}\left( 18200\omega \text{+}16425\right)
f_{n}^{2}\text{+}
\end{equation*}%
\begin{equation}
+\left( 2915\omega +3819\right) f_{n-1}^{2}-\left( 26448\omega +24138\right)
\cdot \left( -1\right) ^{n}].\   \tag{4.24.}
\end{equation}%
Using Remark 4.3 (i,ii) and the definition of the generalized Fibonacci
numbers, it is easy to compute that: 
\begin{equation*}
E\left( f_{n+1},f_{n+4},f_{n+7}\right) \text{+}E\left(
f_{n},f_{n+1},f_{n+2}\right) \text{+}E\left( f_{n+3},f_{n+4},f_{n+5}\right) 
\text{+}E\left( f_{n+6},f_{n+7},f_{n+8}\right)
\end{equation*}%
\begin{equation*}
+E\left( \omega f_{n},f_{n+5},f_{n+7}\right) \text{+}E\left(
f_{n+1},f_{n+3},\omega f_{n+8}\right) \text{+}E\left( f_{n+2},f_{n+4},\omega
f_{n+6}\right) \text{+}
\end{equation*}%
\begin{equation*}
E\left( f_{n},f_{n+4},\omega ^{2}f_{n+8}\right) \text{+}E\left( \omega
^{2}f_{n+1},f_{n+5},f_{n+6}\right) \text{+}E\left( f_{n+2},f_{n+3},\omega
^{2}f_{n+7}\right) =
\end{equation*}%
\begin{equation*}
=f_{n+2}\cdot \left( \omega h_{2n}^{30766,27923}\text{+}h_{2n}^{22358,20533}%
\right) \text{+}f_{n+3}\cdot \left( \omega h_{2n}^{4368,1453}\text{+}%
h_{2n}^{14128,12163}\right) -
\end{equation*}%
\begin{equation*}
\text{-}f_{n}^{2}\cdot \left( \omega h_{n+3}^{45013,22563}\text{+}%
h_{n+3}^{33683,27523}\right) +\left( -1\right) ^{n+1}\cdot \left( \omega
h_{n+3}^{1472,26448}+h_{n+3}^{12982,24138}\right) .
\end{equation*}%
Therefore, we obtain: 
\begin{equation*}
E\left( f_{n+1},f_{n+4},f_{n+7}\right) \text{+}E\left(
f_{n},f_{n+1},f_{n+2}\right) \text{+}E\left( f_{n+3},f_{n+4},f_{n+5}\right) 
\text{+}E\left( f_{n+6},f_{n+7},f_{n+8}\right) \text{+}
\end{equation*}%
\begin{equation*}
E\left( \omega f_{n},f_{n+5},f_{n+7}\right) +E\left( f_{n+1},f_{n+3},\omega
f_{n+8}\right) +E\left( f_{n+2},f_{n+4},\omega f_{n+6}\right) +
\end{equation*}%
\begin{equation*}
E\left( f_{n},f_{n+4},\omega ^{2}f_{n+8}\right) +E\left( \omega
^{2}f_{n+1},f_{n+5},f_{n+6}\right) +E\left( f_{n+2},f_{n+3},\omega
^{2}f_{n+7}\right) =
\end{equation*}%
\begin{equation*}
f_{n+2}\cdot \left( \omega h_{2n}^{30766,27923}\text{+}h_{2n}^{22358,20533}%
\right) \text{+}f_{n+3}\cdot \left( \omega h_{2n}^{4368,1453}\text{+}%
h_{2n}^{14128,12163}\right) \text{-}
\end{equation*}%
\begin{equation}
\text{-}f_{n}^{2}\cdot \left( \omega h_{n+3}^{45013,22563}\text{+}%
h_{n+3}^{33683,27523}\right) \text{+}\left( \text{-}1\right) ^{n+1}\cdot
\left( \omega h_{n+3}^{1472,26448}\text{+}h_{n+3}^{12982,24138}\right) .\  
\tag{4.25.}
\end{equation}%
Adding equalities $(4.7),(4.15),(4.25),$ we have: 
\begin{equation*}
\eta (F_{n})\text{=}4a^{2}h_{n+3}^{211,14}\cdot \left( h_{2n}^{84,135}\text{-%
}2f_{n}^{2}\right) \text{+}8h_{n+3}^{8,3}\cdot \left( h_{2n}^{12,20}\text{-}%
f_{n}^{2}\right) \text{+}
\end{equation*}%
\begin{equation*}
+a[f_{n+2}\left( \omega h_{2n}^{30766,27923}\text{+}h_{2n}^{22358,20533}%
\right) \text{+}
\end{equation*}%
\begin{equation*}
f_{n+3}\cdot \left( \omega h_{2n}^{4368,1453}\text{+}h_{2n}^{14128,12163}%
\right) \text{-}f_{n}^{2}\left( \omega h_{n+3}^{45013,22563}\text{+}%
h_{n+3}^{33683,27523}\right) \text{-}
\end{equation*}%
\begin{equation*}
-\left( \text{-}1\right) ^{n}\left( \omega h_{n+3}^{1472,26448}\text{+}%
h_{n+3}^{12982,24138}\right) ].
\end{equation*}%
$\Box \medskip $\smallskip \newline
\textbf{Corollary 4.6. } \textit{Let} $F_{n}$ \textit{be the n th Fibonacci
symbol element. Let }$\omega $\textit{\ be a primitive root of order }$3$%
\textit{\ of unity and let }$K=\mathbb{Q}\left( \omega \right) $ \textit{be} 
\textit{the cyclotomic field $.$ Then, the norm of} $F_{n}$ \textit{is} 
\begin{equation*}
\eta (F_{n})\text{=}f_{n+2}\left( \omega h_{2n}^{30766,27923}\text{+}%
h_{2n}^{26822,27753}\right) \text{+}f_{n+3}\cdot \left( \omega
h_{2n}^{4368,1453}\text{+}h_{2n}^{19120,20203}\right) \text{-}
\end{equation*}%
\begin{equation*}
\text{-}f_{n}^{2}\cdot (\omega h_{n+3}^{45013,22563}\text{+}%
h_{n+3}^{33835,27659})\text{+}\left( \text{-}1\right) ^{n+1}\left( \omega
h_{n+3}^{1472,26448}\text{+}h_{n+3}^{12982,24138}\right) .
\end{equation*}%
\smallskip \newline
\textbf{Proof.} From the relation $(4.7),$ we know that 
\begin{equation*}
E\left( f_{n+2},f_{n+5},f_{n+8}\right) =4\left( 11f_{n+2}+14f_{n+3}\right)
\cdot \left( h_{2n}^{84,135}-2f_{n}^{2}\right) .
\end{equation*}%
Therefore, we obtain 
\begin{equation}
E\left( f_{n+2},f_{n+5},f_{n+8}\right) =f_{n+2}\left(
h_{2n}^{3696,5940}-88f_{n}^{2}\right) +f_{n+3}\left(
h_{2n}^{4704,7560}-112f_{n}^{2}\right) .\   \tag{4.26.}
\end{equation}%
From the relation $(4.15),$ we know that 
\begin{equation*}
E\left( f_{n},f_{n+3},f_{n+6}\right) =8\left( 8f_{n+2}+3f_{n+3}\right) \cdot
\left( h_{2n}^{12,20}-f_{n}^{2}\right) .
\end{equation*}%
Then, we obtain 
\begin{equation}
E\left( f_{n},f_{n+3},f_{n+6}\right) =f_{n+2}\left(
h_{2n}^{768,1280}-64f_{n}^{2}\right) +f_{n+3}\left(
h_{2n}^{288,480}-24f_{n}^{2}\right) .\   \tag{4.27.}
\end{equation}%
Adding equalities (4.26) and (4.27) it is easy to compute that: 
\begin{equation*}
E\left( f_{n+2},f_{n+5},f_{n+8}\right) \text{+}E\left(
f_{n},f_{n+3},f_{n+6}\right) \text{=}
\end{equation*}%
\begin{equation*}
\text{=}f_{n+2}\left( h_{2n}^{3696,5940}\text{+}h_{2n}^{768,1280}\right) 
\text{+}f_{n+3}\left( h_{2n}^{4704,7560}\text{+}h_{2n}^{288,480}\right) 
\text{-}f_{n}^{2}\left( 152f_{n+2}\text{+}136f_{n+3}\right) \text{=}
\end{equation*}%
\begin{equation*}
\text{=}f_{n+2}h_{2n}^{4464,7220}\text{+}f_{n+3}h_{2n}^{4992,8040}\text{-}%
f_{n}^{2}h_{n+3}^{152,136}.
\end{equation*}%
It results 
\begin{equation*}
\eta (F_{n})\text{=}f_{n+2}h_{2n}^{4464,7220}\text{+}f_{n+3}h^{4992,8040}%
\text{-}f_{n}^{2}h_{n+3}^{152,136}\text{+}
\end{equation*}%
\begin{equation*}
\text{+}f_{n+2}\left( \omega h_{2n}^{30766,27923}\text{+}%
h_{2n}^{22358,20533}\right) \text{+}f_{n+3}\cdot
\end{equation*}%
\begin{equation*}
\left( \omega h_{2n}^{4368,1453}\text{+}h_{2n}^{14128,12163}\right) \text{-}%
f_{n}^{2}\left( \omega h_{n+3}^{45013,22563}\text{+}h_{n+3}^{33683,27523}%
\right) \text{-}
\end{equation*}%
\begin{equation*}
-\left( \text{-}1\right) ^{n}\left( \omega h_{n+3}^{1472,26448}\text{+}%
h_{n+3}^{12982,24138}\right) .
\end{equation*}%
Therefore%
\begin{equation*}
\eta (F_{n})\text{=}f_{n+2}\left( \omega h_{2n}^{30766,27923}\text{+}%
h_{2n}^{26822,27753}\right) \text{+}f_{n+3}\cdot \left( \omega
h_{2n}^{4368,1453}\text{+}h_{2n}^{19120,20203}\right) -
\end{equation*}%
\begin{equation*}
\text{-}f_{n}^{2}\cdot (\omega h_{n+3}^{45013,22563}\text{+}%
h_{n+3}^{33835,27659})\text{+}\left( \text{-}1\right) ^{n+1}\left( \omega
h_{n+3}^{1472,26448}\text{+}h_{n+3}^{12982,24138}\right) .\newline
\end{equation*}%
\smallskip \newline
In conclusion, even indices of the top of generalized Fibonacci numbers are
very large, the expressions of the norm $\eta (F_{n})$ from Proposition 4.5
and from Corollary 4.6 are much shorten than the formula founded in [Pi; 82]
and, in addition, the powers of Fibonacci numbers in these expressions are $%
1 $ or $2.\Box \medskip $\smallskip \newline

Let $S=\left( \frac{1,1}{\mathbb{Q},\omega }\right) $ be a symbol algebra of
degree $3.$ Using the norm form given in the Corollary 4.6, we obtain:%
\newline
$\eta \left( F_{n}\right) =$$\omega \left(
f_{n+2}h_{2n}^{30766,27923}+f_{n+3}h_{2n}^{4368,1453}-f_{n}^{2}h_{n+3}^{45013,22563}+\left( -1\right) ^{n+1}h_{n+3}^{1472,26448}\right) 
$\newline
$+\left(
f_{n+2}h_{2n}^{26822,27753}+f_{n+3}h_{2n}^{19120,20203}-f_{n}^{2}h_{n+3}^{33835,27659}+\left( -1\right) ^{n+1}h_{n+3}^{12982,24138}\right) . 
$ \newline
If $\eta \left( F_{n}\right) \neq 0,$ we know that the element $F_{n}$ is
invertible. From relation $\eta \left( F_{n}\right) \neq 0,$ it results 
\begin{equation}
~f_{n+2}h_{2n}^{30766,27923}\text{+}f_{n+3}h_{2n}^{4368,1453}\text{-}%
f_{n}^{2}h_{n+3}^{45013,22563}\text{+}\left( -1\right)
^{n+1}h_{n+3}^{1472,26448}\neq 0  \tag{4.28.}
\end{equation}
or 
\begin{equation}
f_{n+2}h_{2n}^{26822,27753}\text{+}f_{n+3}h_{2n}^{19120,20203}\text{-}%
f_{n}^{2}h_{n+3}^{33835,27659}\text{+}\left( -1\right)
^{n+1}h_{n+3}^{12982,24138}\neq 0  \tag{4.29.}
\end{equation}
Using relation $\left( 4.1\right) $ and Remark 4.3, i) and v), we obtain that%
\newline
$%
f_{n+2}h_{2n}^{26822,27753}+f_{n+3}h_{2n}^{19120,20203}-f_{n}^{2}h_{n+3}^{33835,27659}+\left( -1\right) ^{n+1}h_{n+3}^{12982,24138}= 
$\newline
$%
=26822f_{n+2}f_{2n-1}+27753f_{n+2}f_{2n}+19120f_{n+3}f_{2n-1}+20203f_{n+3}f_{2n}- 
$\newline
$-33835f_{n}^{2}f_{n+2}-27659f_{n}^{2}f_{n+3}+\left( -1\right)
^{n+1}12982f_{n+2}+\left( -1\right) ^{n+1}24138f_{n+3}=$\newline
$=26822f_{n+2}f_{n-1}^{2}+26822f_{n+2}f_{n}^{2}+27753f_{n+2}f_{n}^{2}+2\cdot
27753f_{n-1}f_{n}f_{n+2}+$\newline
$+19120f_{n+3}f_{n}^{2}+19120f_{n+3}f_{n-1}^{2}+20203f_{n+3}f_{2n}-$\newline
$-33835f_{n}^{2}f_{n+2}-27659f_{n}^{2}f_{n+3}+\left( -1\right)
^{n+1}12982f_{n+2}+\left( -1\right) ^{n+1}24138f_{n+3}=$\newline
$=26822f_{n+2}f_{n-1}^{2}+\left( -1\right) ^{n+1}12982f_{n+2}+
26822f_{n+2}f_{n}^{2}+27753f_{n+2}f_{n}^{2}-$ \newline
$-27659f_{n}^{2}f_{n+3}+2\cdot 27753f_{n-1}f_{n}f_{n+2}+\left( -1\right)
^{n+1}24138f_{n+3}+$\newline
$%
+19120f_{n+3}f_{n}^{2}+20203f_{n+3}f_{2n}-33835f_{n}^{2}f_{n+2}+19120f_{n+3}f_{n-1}^{2}. 
$\newline
We remark that $26822f_{n+2}f_{n-1}^{2}+\left( -1\right)
^{n+1}12982f_{n+2}>0.$\newline
We have $26822f_{n+2}f_{n}^{2}+27753f_{n+2}f_{n}^{2}-27659f_{n}^{2}f_{n+3}=$%
\newline
$=54575f_{n+2}f_{n}^{2}-27659f_{n}^{2}f_{n+3}\geq 0,$ since $%
54575f_{n+2}>27659f_{n+3}$ is equivalently with $\frac{f_{n+2}}{f_{n+3}}>%
\frac{27659}{54575}\approx 0,506,$which is true for all $n\geq 1.$For $n=0,$
we have $54575f_{n+2}f_{n}^{2}-27659f_{n}^{2}f_{n+3}=0.$\newline
Obviously, $2\cdot 27753f_{n-1}f_{n}f_{n+2}+\left( -1\right)
^{n+1}24138f_{n+3}=$\newline
$=55506f_{n-1}f_{n}f_{n+2}+\left( -1\right) ^{n+1}24138f_{n+3}>0,$ since $%
\frac{f_{n+3}}{f_{n+2}}<\frac{55506}{24138}\approx 2,29.$\newline
And, finally, $%
19120f_{n+3}f_{n}^{2}+20203f_{n+3}f_{2n}-33835f_{n}^{2}f_{n+2}>0,$ since $%
f_{2n}>f_{n}^{2}$ and $19120+20203=39323>33835.$ \newline
From the above, we have that the term\newline
$%
f_{n+2}h_{2n}^{26822,27753}+f_{n+3}h_{2n}^{19120,20203}-f_{n}^{2}h_{n+3}^{33835,27659}+\left( -1\right) ^{n+1}h_{n+3}^{12982,24138}>0 
$ for all $n\in \mathbb{N}.$

We just proved the following result:\medskip

\textbf{Theorem 4.7. }\ \textit{Let} \ $S=\left( \frac{1,1}{\mathbb{Q}%
,\omega }\right) $ \textit{be a symbol algebra of degree} $3.$ \textit{All
symbol Fibonacci elements are invertible elements.~}$\Box \medskip $%
\smallskip \smallskip \newline
\begin{equation*}
\end{equation*}%
\textbf{Conclusions.} In this paper, we studied some properties of the
matrix representation of symbol algebras of degree $3.$ Using some
properties of left matrix representations, \ we solved equations with
coefficients in these algebras. The study of symbol algebras of degree three
involves very complicated calculus and, usually, can be hard to find \
examples for some notions. We introduced the Fibonacci symbol elements , we
gave an easier expression of reduced norm of Fibonacci symbol elements, and,
from this formula, we found examples of many invertible elements, namely all
Fibonacci symbol elements are invertible. Using some ideas from this paper,
we expect to obtain some interesting new results in a further research and
we hope that this kind of sets of invertible elements will help us to
provide, in the future, examples of symbol division algebra of degree three
or of degree greater than three. \newline
We think that this thing is possible, since in the case of quaternion
algebras the set obtained by union of the set of Fibonacci quaternion
elements with the set of Lucas quaternion elements and with the set of
Mersenne quaternion elements is very close to being an algebra. This idea
will be developped in our next researchs.%
\begin{equation*}
\end{equation*}%
\textbf{Acknowledgments.} The authors would like to thank the referee for
his/her patience and suggestions, which helped us improve this paper.\qquad
\smallskip \newline
\begin{equation*}
\end{equation*}%
\textbf{References}

\begin{equation*}
\end{equation*}

[Fa; 88]$~$\ J.R. Faulkner, \textit{Finding octonion algebras in associative
algebras,} Proc. Amer. Math. Soc. \textbf{104}(4)(1988), 1027-1030.

[Fla; 12] R. Flatley, \textit{Trace forms of Symbol Algebras}, Algebra
Colloquium, \textbf{19}(2012), 1117-1124.

[Fl; Sh; 12] C. Flaut, V. Shpakivskyi, \textit{On Generalized Fibonacci
Quaternions and Fibonacci-Narayana Quaternions,} Adv. Appl. Clifford
Algebras, \textbf{23(3)}(2013), 673-688.

[Gi, Sz; 06 ] P. Gille, T. Szamuely, \textit{Central Simple Algebras and
Galois Cohomology}, Cambridge University Press, 2006.

[Ha; 00] M. Hazewinkel, \textit{Handbook of Algebra}, Vol. 2, North Holland,
Amsterdam, 2000.

[Ho; 61] A. F. Horadam, \textit{A Generalized Fibonacci Sequence}, Amer.
Math. Monthly, \textbf{68}(1961), 455-459.

[Ho; 63] A. F. Horadam, \textit{Complex Fibonacci Numbers and Fibonacci
Quaternions}, Amer. Math. Monthly, \textbf{70}(1963), 289-291.

[Ka; 53] L. Kaplansky, \textit{Infinite-dimensional quadratic forms
admitting composition, }Proc. Amer. Math., \textbf{\ 4(}1953), 956-960.

[La; 05] T. Y. Lam, \textit{Introduction To Quadratic Forms Over Fields},
Graduate Studies in Mathematics, vol 67, American Mathematical Society,
Providence RI, 2005.

[Mil; 12] J.S. Milne,\textit{\ Class field theory}
(http://www.math.lsa.umich.edu/jmilne);

[Mi; 71] \ J. Milnor, \ \textit{Introduction to Algebraic K-Theory, }Annals
of Mathematics Studies, Princeton Univ. Press, 1971.

[Pi; 82] R.S, Pierce, \textit{Associative Algebras}, Springer Verlag, New
York, Heidelberg, Berlin, 1982.

[Sa; Fa; Ci; 09] D. Savin, C. Flaut, C. Ciobanu, \textit{Some properties of
the symbol algebras}, Carpathian Journal of Mathematics , vol. 25, No. 2
(2009), 239-245.

[Ti; 00] Y. Tian, \textit{Matrix reprezentations of octonions and their
applications, }Adv. in Appl. Clifford Algebras, \textbf{10}(1)( 2000),
61-90. \bigskip 
\begin{equation*}
\end{equation*}

Cristina FLAUT

{\small Faculty of Mathematics and Computer Science,}

{\small Ovidius University,}

{\small Bd. Mamaia 124, 900527, CONSTANTA, ROMANIA}

{\small http://cristinaflaut.wikispaces.com/}

{\small http://www.univ-ovidius.ro/math/}

{\small e-mail: cflaut@univ-ovidius.ro, cristina\_flaut@yahoo.com}

\bigskip

Diana SAVIN

{\small Faculty of Mathematics and Computer Science, }

{\small Ovidius University, }

{\small Bd. Mamaia 124, 900527, CONSTANTA, ROMANIA }

{\small http://www.univ-ovidius.ro/math/ }

{\small e-mail: \ savin.diana@univ-ovidius.ro, \ dianet72@yahoo.com}

\end{document}